# From the Half-Order Recurrence to General Fractional-Order Differentiation on Monomials: Functional Continuation and Operator Composition

Davit Kapanadze
Georgian National University SEU, Tbilisi, Georgia
Associate Member, Institute of Mathematics and its Applications (IMA), United Kingdom
ORCID: 0009-0006-8112-5084
d.kapanadze@seu.edu.ge
DavitKapanadze06@gmail.com

**Abstract**

The present work develops a general-order fractional differentiation operator on monomials from the discrete half-order coefficient structure established previously in Part I. The results of the half-order construction are taken here as the starting point and are not rederived.

The first objective is to pass from the discrete coefficient law to a functional relation defined on a broader domain. The one-step transfer law obtained from the integer coefficient family is extended from integer arguments to a real argument, and a natural canonical Gamma-functional continuation is then selected under the stated normalization and regularity requirements. In this construction, the Gamma function is not introduced as the initial definition of fractional differentiation; rather, it appears as the functional form associated with the coefficient law inherited from the discrete construction.

The second objective is to extend the operator order beyond one half. The first derivative is divided into equal operator steps, beginning with specific rational orders and proceeding to orders $1/m$, $k/m$, and subsequently to general rational and real orders. The resulting coefficient is independently verified through the addition law for operator orders on the monomial system, subject to the requirement that all intermediate expressions occurring in the composition are defined.

On the corresponding parameter domain, the construction yields the standard Gamma-ratio monomial expression for fractional differentiation. For nonnegative integer orders it reduces to ordinary repeated differentiation, while for positive non-integer orders it agrees at the monomial level with the left-sided Riemann–Liouville formula with lower limit 0. The relation with the Caputo operator is also discussed, including the distinction that remains for constants and low-degree polynomial terms.

The purpose of the work is not to introduce a new classical fractional derivative or to replace the Riemann–Liouville or Caputo theories. Its purpose is to provide a constructive algebraic-operator route from a discrete coefficient law to the general fractional-order monomial formula, making explicit the roles of functional continuation, Gamma normalization, and operator composition.

**Keywords:** fractional calculus; general-order fractional derivative; monomials; algebraic-operator construction; functional continuation; Gamma function; operator composition; addition law for operators.

## Introduction

The ordinary differentiation operator acts on a monomial by the rule

$$Dx^{\beta} = \beta x^{\beta-1}.$$

Under repeated differentiation of integer order, the operators obey the addition law

$$D^m D^n = D^{m+n}$$

on the corresponding functional domain. This algebraic property raises a natural question: whether the order of differentiation can be extended beyond the integers while preserving an appropriate composition structure and, at the same time, whether the coefficients appearing in such an extension can be obtained constructively rather than introduced in their final form from the outset.

The classical theory of fractional differentiation – including the Riemann–Liouville and Caputo operators, their integral representations, and their applications – is developed in standard references [1–4, 10, 35], discussed historically in [14, 18], and treated in a wide range of applied settings [5–13, 25]. More general fractional integrals and derivatives based on convolution structures and Sonine kernels have also been developed; in particular, the works of Y. Luchko [38–40] provide a broad framework for such operators and their operational calculus.

The Gamma-functional formula for the fractional derivative of a monomial belongs to the classical historical development of fractional calculus, going back to Euler and later developments associated with Liouville and other nineteenth-century constructions [14, 18]. The present work therefore does not claim a new final monomial formula. Its aim is to reach the known Gamma-functional expression by a different, constructive algebraic-operator route.

The starting point of the present investigation is the half-order algebraic-operator construction established in Part I [26] (arXiv:2607.19482), where the coefficient recurrence, the integer and half-integer coefficient chains, the appearance of double factorials, and the normalization were obtained from the half-step composition requirement. These results are used here as established facts and are not rederived. Part I, in turn, continues the author's earlier limit-free algebraic-geometric approach to differentiation [15–17].

The present paper begins where Part I ends. The first question is whether the discrete coefficient law obtained there admits a natural functional continuation from its original discrete argument set to a broader real domain. Instead of reconstructing the discrete coefficient families, we extract from them the one-step transfer relation required for this continuation. This relation leads to a functional equation whose simplest canonical realization is expressed through the Gamma function.

The second question concerns the order of the operator itself. The half-order case represents a division of the first derivative into two equal operator steps. The natural next problem is to investigate analogous divisions into three, four, or, in general, $m$ equal steps and to determine whether these constructions can be organized into operators of rational and subsequently real order. The addition law for operator orders then provides an independent algebraic mechanism for determining and verifying the general coefficient.

The scope of the construction is intentionally limited. Agreement with the Riemann–Liouville formula on monomials does not imply that the complete Riemann–Liouville theory has been reconstructed. No general integral or convolution kernel, general mechanism of non-locality, broad function-space theory, or complete theory of fractional differential equations is constructed here. The operator family developed in this paper acts first on monomials and then, by linearity, on their finite linear combinations.

The paper is organized as follows. Chapter I develops the functional continuation of the coefficient law inherited from Part I and its Gamma-functional representation. Chapter II motivates the transition from the half-order construction to operators of more general order. Chapter III obtains and verifies the general coefficient through operator composition and the addition law. Chapter IV considers the linear extension to polynomial initial functions and comparison with the Riemann–Liouville and Caputo monomial rules. Chapter V provides computational and visual examples, and Chapter VI summarizes the results, limitations, and directions for further development.

## Convention of Parameters and Notation

Throughout the paper the following conventions are used, unless a different condition is explicitly stated at a particular point.

**1. The variable $x$.** We take

$$x > 0,$$

so that real non-integer powers $x^{\beta}$ and $x^{\beta-\alpha}$ are real and single-valued. The case $x \leq 0$ is not considered in the present construction.

**2. The monomial exponent $\beta$.** The symbol $\beta$ denotes the exponent of the original monomial $x^{\beta}$. The basic monomial domain on which the classical left-sided Riemann–Liouville formula with lower limit 0 applies directly is

$$\beta > -1.$$

When $\beta+1-\alpha$ is a non-positive integer, the coefficient is interpreted, where appropriate, through the reciprocal Gamma function. Cases in which $\Gamma(\beta+1)$ itself has a pole are not included in the basic monomial domain and require a separate meromorphic or distributional treatment.

**3. Operator orders $\alpha$ and $\gamma$.** These symbols denote operator orders. Positive values are treated as differentiation orders, while a negative order is interpreted, at the monomial level, as a Riemann–Liouville-type integral of the corresponding positive order. At the constructive stage for positive rational orders we use

$$\left(D^{1/m}\right)^m = D,$$

where the first derivative is divided into $m$ equal operator steps.

**4. Integer indices $n, k, m, j$.** We use

$$n, j, k \in \mathbb{N}_0, \qquad m \in \mathbb{N}, \qquad m \geq 2.$$

Here $n$ denotes a nonnegative integer monomial degree, $m$ is the number of equal operator steps and the denominator of the rational order $k/m$, $k$ is the number of successive steps and the numerator of $k/m$, and $j$ is the index of a polynomial term.

**5. The addition law for operators.** The identity

$$D^{\alpha} D^{\gamma} x^{\beta} = D^{\alpha+\gamma} x^{\beta}$$

is used only for those values of $\beta$, $\alpha$, and $\gamma$ for which all intermediate monomials and all corresponding coefficients participating in the composition are defined. This restriction is substantive and is understood whenever the addition law is invoked.

**6. Poles of the Gamma function and the reciprocal Gamma function.** If

$$\beta + 1 - \alpha \in \{0, -1, -2, \dots\},$$

then, where needed, the coefficient is interpreted through the reciprocal Gamma function,

$$\frac{1}{\Gamma(\beta + 1 - \alpha)} = 0.$$

This convention does not mean that $\Gamma(0)$, $\Gamma(-1)$, and similar values are finite; it only specifies the value of the reciprocal-Gamma coefficient at those points.

**7. Behavior as $\mathbf{x \to 0^+}$.** The general monomial expression contains the power $\beta - \alpha$. Hence the result is regular as $x \to 0^+$ when $\beta - \alpha > 0$, constant when $\beta - \alpha = 0$, and generally singular when $\beta - \alpha < 0$, except in cases where the Gamma coefficient vanishes through the reciprocal-Gamma convention. Graphical and numerical examples are therefore taken on strictly positive x whenever necessary.

# Chapter I – Functional Continuation: From $c(\beta)$ to the Gamma Function

## 1.1. Starting Point from Part I and the One-Step Transfer Law

The discrete half-order coefficient structure required in this chapter was established in Part I [26]. In particular, Part I determined the normalized coefficient families on the integer and half-integer sectors of the graded monomial space. Since their derivation is already available there, it will not be repeated in the present work.

For the present purpose we need only one consequence of the integer coefficient family. Taking the ratio of two successive normalized integer coefficients obtained in Part I gives

$$c(n+1) = \frac{n+1}{n+\frac{1}{2}} c(n), \qquad n = 0, 1, 2, \dots \tag{1.1}$$

This one-step relation is the only discrete result from Part I required for the functional continuation developed below. The problem is now to determine how this relation may be extended from nonnegative integer arguments to a continuous argument.

## 1.2. Formation of the Functional Equation

Formula (1.1) was originally obtained only for nonnegative integer $n$. However, its right-hand side contains only ordinary algebraic operations:

$$\frac{n+1}{n+\frac{1}{2}}\, c\,(n)\,.$$

This expression retains meaning not only for integer but also for the corresponding real arguments. Hence, as the requirement of functional continuation, we take that the unknown function $c\,(\beta)$ should preserve the same transfer law:

$$c\,(\beta+1) = \frac{\beta+1}{\beta+\frac{1}{2}}\, c\,(\beta)\,. \tag{1.2}$$

Formula (1.2) is used for those $\beta$ for which the expression is defined. This is the basic functional equation. It states that when the argument is increased by one, $c\,(\beta)$ is multiplied by the factor

$$\frac{\beta+1}{\beta+\frac{1}{2}}.$$

Formula (1.2) should be understood not as the automatic consequence of the unique possible continuation from the discrete data, but as a natural principle of continuation: on the real argument as well, we preserve the same one-step algebraic law that has already been proved on integer values.

Now the task is to find such a natural solution of this functional equation that agrees with the normalization

$$c\,(0) = \frac{1}{\sqrt{\pi}}. \tag{1.3}$$

## 1.3. The Gamma Function and Its Transfer Law

The functional equation (1.2) contains two linear factors, $\beta+1$ and $\beta+\frac{1}{2}$. To represent these factors naturally, we use the Gamma function. The Gamma function may be regarded as the canonical continuation of the factorial law to non-integer arguments. It is normalized by the condition $\Gamma\,(1) = 1$ and satisfies the transfer law

$$\Gamma\,(z+1) = z\Gamma\,(z)\,. \tag{1.4}$$

Here $\Gamma\,(1) = 1$ is the initial normalization condition. It is precisely on the basis of this condition and the transfer law that the Gamma function coincides on positive integer arguments with the ordinary factorial. Indeed,

$$\Gamma\,(2) = 1\Gamma\,(1) = 1, \qquad \Gamma\,(3) = 2\Gamma\,(2) = 2!, \qquad \Gamma\,(4) = 3\Gamma\,(3) = 3!,$$

and, in general, for every nonnegative integer $n$,

$$\Gamma\,(n+1) = n!. \tag{1.5}$$

Thus, the Gamma function preserves the basic recurrent law of the factorial and extends it to non-integer arguments as well. For example, it gives meaning to expressions such as $\Gamma\left(\frac{1}{2}\right)$, $\Gamma\left(\frac{3}{2}\right)$, and $\Gamma\left(\beta+\frac{1}{2}\right)$.

It should be noted that the functional equation $F(z+1) = zF(z)$ alone is not sufficient to determine a unique function. Additional conditions are needed to fix it. On the positive real axis, the normalized Gamma function is characterized by the transfer law together with positivity and logarithmic convexity; this is the classical canonical continuation of the factorial [30].

Substituting $z = \beta + 1$ into the transfer law, we obtain

$$\Gamma(\beta+2) = (\beta+1)\,\Gamma(\beta+1), \tag{1.6}$$

from which

$$\beta + 1 = \frac{\Gamma(\beta+2)}{\Gamma(\beta+1)}.$$

Analogously, for $z = \beta + \frac{1}{2}$,

$$\Gamma\left(\beta+\frac{3}{2}\right) = \left(\beta+\frac{1}{2}\right)\Gamma\left(\beta+\frac{1}{2}\right), \tag{1.7}$$

and therefore

$$\beta + \frac{1}{2} = \frac{\Gamma\left(\beta+\frac{3}{2}\right)}{\Gamma\left(\beta+\frac{1}{2}\right)}.$$

Thus, both linear factors present in the functional equation are expressed by ratios of neighboring Gamma-function values. It is precisely this circumstance that makes the Gamma function a natural instrument for solving the given functional equation.

### 1.4. Construction of the Gamma-Functional Candidate

The form of the functional equation

$$c(\beta+1) = \frac{\beta+1}{\beta+\frac{1}{2}}\,c(\beta)$$

suggests that we seek $c(\beta)$ in the form of a ratio of two Gamma functions. Consider the candidate

$$c(\beta) = \frac{\Gamma(\beta+1)}{\Gamma\left(\beta+\frac{1}{2}\right)}. \tag{1.8}$$

Let us check whether it satisfies the functional equation. Substituting the argument $\beta + 1$, we obtain

$$c(\beta+1) = \frac{\Gamma(\beta+2)}{\Gamma\left(\beta+\frac{3}{2}\right)}. \tag{1.9}$$

Using (1.6) and (1.7),

$$\Gamma(\beta+2) = (\beta+1)\,\Gamma(\beta+1),$$

and

$$\Gamma\left(\beta+\frac{3}{2}\right)=\left(\beta+\frac{1}{2}\right)\Gamma\left(\beta+\frac{1}{2}\right).$$

Hence

$$c\left(\beta+1\right)=\frac{\beta+1}{\beta+\frac{1}{2}}\frac{\Gamma\left(\beta+1\right)}{\Gamma\left(\beta+\frac{1}{2}\right)}=\frac{\beta+1}{\beta+\frac{1}{2}}\,c\left(\beta\right),$$

which exactly coincides with the functional equation (1.2). Accordingly, (1.8) is a natural Gamma-functional solution of the functional equation.

We call this solution the natural, or canonical, Gamma-functional continuation – not the unique solution of the functional equation. The precise meaning of this choice and the question of uniqueness are clarified in Section 1.10.

## 1.5. Verification of the Formula on Integer Arguments

We must verify that the Gamma-functional formula exactly recovers, for integer $n$, the coefficients obtained in Part I. Let

$$\beta=n,\qquad\qquad n=0,1,2,\ldots.$$

Then

$$c\left(n\right)=\frac{\Gamma\left(n+1\right)}{\Gamma\left(n+\frac{1}{2}\right)}.\tag{1.10}$$

For the numerator,

$$\Gamma\left(n+1\right)=n!.\tag{1.11}$$

For the half-integer argument, repeated use of the Gamma transfer law gives

$$\Gamma\left(n+\frac{1}{2}\right)=\frac{(2n-1)!!}{2^n}\Gamma\left(\frac{1}{2}\right).\tag{1.12}$$

Using the classical value $\Gamma\left(\frac{1}{2}\right)=\sqrt{\pi}$, whose agreement with the normalization chosen in Part I will be checked separately in Section 1.7, we obtain

$$\Gamma\left(n+\frac{1}{2}\right)=\frac{(2n-1)!!}{2^n}\sqrt{\pi}.\tag{1.13}$$

Substituting (1.11) and (1.13) into (1.10),

$$c\left(n\right)=\frac{n!}{\frac{(2n-1)!!}{2^n}\sqrt{\pi}}=\frac{2^n n!}{(2n-1)!!\sqrt{\pi}}.$$

Since $(2n)!!=2^n n!$, we obtain

$$c\,(n) = \frac{(2n)!!}{(2n-1)!!\sqrt{\pi}}.$$

Thus, the Gamma-functional formula on integer arguments exactly recovers the double-factorial coefficients obtained in Part I.

### 1.6. Verification on Half-Integer Arguments

Now let

$$\beta = n - \frac{1}{2}, \qquad n = 1, 2, 3, \ldots.$$

Formula (1.8) gives

$$c\left(n - \frac{1}{2}\right) = \frac{\Gamma\left(n + \frac{1}{2}\right)}{\Gamma\,(n)}. \tag{1.14}$$

Since $\Gamma\,(n) = (n-1)!$ and

$$\Gamma\left(n + \frac{1}{2}\right) = \frac{(2n-1)!!}{2^n}\sqrt{\pi},$$

we have

$$c\left(n - \frac{1}{2}\right) = \frac{(2n-1)!!\sqrt{\pi}}{2^n\,(n-1)!}.$$

Using $(2n-2)!! = 2^{n-1}\,(n-1)!$, hence $2^n\,(n-1)! = 2\,(2n-2)!!$, we obtain

$$c\left(n - \frac{1}{2}\right) = \frac{(2n-1)!!\sqrt{\pi}}{2\,(2n-2)!!}. \tag{1.15}$$

This coincides exactly with the normalized half-integer coefficient obtained in Part I. For the first two values,

$$c\left(\frac{1}{2}\right) = \frac{\sqrt{\pi}}{2}, \qquad c\left(\frac{3}{2}\right) = \frac{3\sqrt{\pi}}{4}.$$

Thus, same Gamma-functional formula describes uniformly both the integer and the half-integer sectors.

### 1.7. Agreement with the Normalization and $\Gamma\,(1/2) = \sqrt{\pi}$

In formula (1.8), let $\beta = 0$. We obtain

$$c\,(0) = \frac{\Gamma\,(1)}{\Gamma\left(\frac{1}{2}\right)}.$$

Since $\Gamma(1) = 1$,

$$c(0) = \frac{1}{\Gamma\left(\frac{1}{2}\right)}. \tag{1.16}$$

In Part I, on the basis of compatibility with the Riemann–Liouville monomial formula, the normalization

$$c(0) = \frac{1}{\sqrt{\pi}} \tag{1.17}$$

was chosen.

Comparing (1.16) and (1.17),

$$\frac{1}{\Gamma\left(\frac{1}{2}\right)} = \frac{1}{\sqrt{\pi}},$$

from which

$$\Gamma\left(\frac{1}{2}\right) = \sqrt{\pi}. \tag{1.18}$$

It is important to interpret this step correctly. Formula (1.18), in this text, does not constitute an independent proof of the classical identity $\Gamma(1/2) = \sqrt{\pi}$. The normalization $c(0) = 1/\sqrt{\pi}$ was chosen in Part I for compatibility with the classical Riemann–Liouville monomial formula. Hence, (1.18) expresses the agreement of the Gamma-functional continuation with the already-chosen normalization. The logical sequence is therefore:

$$\text{algebraic recurrence} \to \text{chosen normalization} \to c(0) = \frac{1}{\sqrt{\pi}} \to c(0) = \frac{\Gamma(1)}{\Gamma(1/2)} \to \Gamma\left(\frac{1}{2}\right) = \sqrt{\pi}.$$

The classical evaluation of this constant itself traces back to Wallis's infinite product [19].

## 1.8. The Complete Monomial Formula for the Half-Derivative

On the basis of the functional continuation obtained, the half-order operator on a monomial can be written as

$$D^{1/2}x^{\beta} = c(\beta)\, x^{\beta-1/2},$$

where

$$c(\beta) = \frac{\Gamma(\beta+1)}{\Gamma\left(\beta+\frac{1}{2}\right)}.$$

Thus,

$$D^{1/2}x^{\beta} = \frac{\Gamma(\beta+1)}{\Gamma\left(\beta+\frac{1}{2}\right)} x^{\beta-1/2}, \qquad x > 0. \tag{1.19}$$

This formula:

- on integer $\beta = n$ recovers

$$D^{1/2}x^{n} = \frac{(2n)!!}{(2n-1)!!\sqrt{\pi}} x^{n-1/2};$$

- on half-integer $\beta = n - \frac{1}{2}$ gives

$$D^{1/2}x^{n-1/2} = \frac{(2n-1)!!\sqrt{\pi}}{2\,(2n-2)!!} x^{n-1};$$

- and, on the domain where the Gamma-functional expression is meaningful, also determines the coefficient for a general real $\beta$.

Formula (1.19) coincides with the left-sided Riemann–Liouville half-derivative monomial formula with lower limit 0. In the present chapter, however, it is obtained not from an integral definition but through the functional continuation of the coefficient law inherited from Part I.

## 1.9. Verification of the Composition Condition on the Continued Form

We now check whether formula (1.19) preserves the original composition requirement. The first action gives

$$D^{1/2}x^{\beta} = \frac{\Gamma(\beta+1)}{\Gamma\left(\beta+\frac{1}{2}\right)} x^{\beta-1/2}.$$

For the second action, the constant coefficient is taken outside the operator:

$$D^{1/2}\left(D^{1/2}x^{\beta}\right) = \frac{\Gamma(\beta+1)}{\Gamma\left(\beta+\frac{1}{2}\right)} D^{1/2}x^{\beta-1/2}.$$

Applying (1.19) to the degree $\beta - \frac{1}{2}$ gives

$$D^{1/2}x^{\beta-1/2} = \frac{\Gamma\left(\beta+\frac{1}{2}\right)}{\Gamma(\beta)} x^{\beta-1}.$$

Therefore,

$$D^{1/2}\left(D^{1/2}x^{\beta}\right) = \frac{\Gamma(\beta+1)}{\Gamma(\beta)} x^{\beta-1}.$$

By the Gamma transfer law, $\Gamma(\beta+1) = \beta\Gamma(\beta)$, so

$$D^{1/2}\left(D^{1/2}x^{\beta}\right) = \beta x^{\beta-1} = Dx^{\beta}. \tag{1.20}$$

Thus, the functionally continued coefficient preserves the original composition requirement on the regular composition domain, where all intermediate monomials and Gamma-functional coefficients are defined. Boundary cases in which an intermediate Gamma factor meets a pole are excluded from this statement and must be checked separately.

### 1.10. On the Canonicity of the Functional Continuation

The functional equation obtained above,

$$c\left(\beta+1\right)=\frac{\beta+1}{\beta+\frac{1}{2}}\,c\left(\beta\right),$$

does not, without additional conditions, determine a unique function by itself. Consider the Gamma-functional solution

$$c_0\left(\beta\right)=\frac{\Gamma\left(\beta+1\right)}{\Gamma\left(\beta+\frac{1}{2}\right)}.$$

If $p\left(\beta\right)$ is any 1-periodic function,

$$p\left(\beta+1\right)=p\left(\beta\right),$$

then

$$c\left(\beta\right)=c_0\left(\beta\right)p\left(\beta\right)=\frac{\Gamma\left(\beta+1\right)}{\Gamma\left(\beta+\frac{1}{2}\right)}p\left(\beta\right)$$

also satisfies the same transfer law. Indeed,

$$c\left(\beta+1\right)=\frac{\Gamma\left(\beta+2\right)}{\Gamma\left(\beta+\frac{3}{2}\right)}p\left(\beta+1\right)=\frac{\beta+1}{\beta+\frac{1}{2}}\,c\left(\beta\right).$$

Thus, the functional equation alone does not rule out an additional 1-periodic factor. At the same time, in order for such a continuation to recover exactly the discrete values obtained in Part I on the corresponding integer and half-integer lattice, the periodic factor must satisfy $p\left(\beta\right)=1$ there. This requirement by itself still does not give absolute uniqueness on the entire continuous domain.

In the present construction we choose the simplest case

$$p\left(\beta\right)\equiv 1,$$

and accordingly regard

$$c\left(\beta\right)=\frac{\Gamma\left(\beta+1\right)}{\Gamma\left(\beta+\frac{1}{2}\right)}$$

as the natural, or canonical, Gamma-functional continuation of the given discrete coefficient family.

This choice is based on the following mutually compatible requirements:

- it recovers exactly the coefficients obtained in Part I on nonnegative integer degrees;
- it recovers, by the same formula, the corresponding half-integer coefficient chain;
- it satisfies the one-step transfer law;
- it is compatible with the normalization $c(0) = 1/\sqrt{\pi}$;
- it preserves the original half-step composition requirement on the regular composition domain;
- at the monomial level, it coincides with the coefficient of the left-sided Riemann–Liouville half-derivative with lower limit 0.

Hence, the Gamma-functional form is not presented here as the absolutely unique solution among all possible solutions of the functional equation. Its canonicity means that it is the simplest continuation compatible with the discrete structure, the normalization, and the composition requirement, without an additional periodic factor.

### 1.11. The Γ-Normalized Monomial Basis and the Half-Step Operator Shift

In the preceding sections, for the half-order operator we obtained

$$D^{1/2}x^{\beta} = \frac{\Gamma(\beta+1)}{\Gamma\left(\beta+\frac{1}{2}\right)}x^{\beta-1/2}.$$

The Gamma-functional coefficient can be simplified by an appropriate normalization of the basis.

**Definition 1.11.1.** On the half-axis $x > 0$, define the Gamma-normalized monomial basis functions by

$$e_{\beta}(x) = \frac{x^{\beta}}{\Gamma(\beta+1)},$$

for those $\beta$ for which the corresponding expressions are defined.

Using linearity and the half-order monomial formula,

$$D^{1/2}e_{\beta}(x) = \frac{1}{\Gamma(\beta+1)}D^{1/2}x^{\beta} = \frac{x^{\beta-1/2}}{\Gamma\left(\beta+\frac{1}{2}\right)}.$$

By the definition of the normalized basis,

$$e_{\beta-1/2}(x) = \frac{x^{\beta-1/2}}{\Gamma\left(\beta+\frac{1}{2}\right)}.$$

Thus,

$$D^{1/2}e_{\beta} = e_{\beta-1/2}.$$

This result shows that in the Gamma-normalized basis the coefficient of the half-derivative disappears, and the action of the operator reduces to an exact half-step shift of the index.

Under a double action,

$$D^{1/2}\left(D^{1/2}e_\beta\right) = D^{1/2}e_{\beta-1/2} = e_{\beta-1}.$$

On the other hand, for the ordinary derivative,

$$De_\beta\left(x\right) = \frac{\beta x^{\beta-1}}{\Gamma\left(\beta+1\right)} = \frac{x^{\beta-1}}{\Gamma\left(\beta\right)} = e_{\beta-1}\left(x\right).$$

Accordingly,

$$D^{1/2}\left(D^{1/2}e_\beta\right) = De_\beta,$$

again only on the domain where all intermediate expressions participating in the composition are defined.

**Remark 1.11.2 (relationship between the bases).** The standard monomial basis $x^\beta$ and the Gamma-normalized basis $e_\beta$ are not alternatives to one another, but two representations serving different purposes. The standard basis is needed to show how the coefficients arise from the discrete recurrence of Part I and how that structure grows into the Gamma-functional form. The Gamma-normalized basis displays the operator content more transparently: instead of a coefficient, only the index shift $\beta \mapsto \beta - \frac{1}{2}$ remains.

Such normalized power functions are naturally related to normalized bases used in the literature on fractional systems and causal operators; in this direction, the works of M. Ortigueira [33–36] are especially relevant. At the present stage, this connection is used only as a structural parallel, not as an already-proved part of a causal or convolutional theory.

It is important to maintain the logical sequence: at this point only the half-order shift law

$$D^{1/2}e_\beta = e_{\beta-1/2}$$

has been obtained. The general-order relation

$$D^\alpha e_\beta = e_{\beta-\alpha}$$

is not used here in advance. It will be established only in Chapter III, after the general monomial operator $D^\alpha$ has been independently obtained and verified. Thus, the Gamma-normalized basis in this chapter serves only to simplify the already-proved half-order structure and does not anticipate the general-order construction.

## 1.12. Chapter Summary

Starting from the discrete coefficient family obtained in Part I, the following steps were carried out in this chapter.

First, the one-step recurrence on nonnegative integer arguments was taken as the discrete starting law:

$$c\left(n+1\right) = \frac{n+1}{n+\frac{1}{2}}\,c\left(n\right).$$

It was then recast as the functional-continuation requirement

$$c(\beta+1) = \frac{\beta+1}{\beta+\frac{1}{2}}\,c(\beta)\,.$$

Using the transfer property of the Gamma function, we obtained the canonical Gamma-functional solution

$$c(\beta) = \frac{\Gamma(\beta+1)}{\Gamma\left(\beta+\frac{1}{2}\right)}.$$

It was confirmed that this formula recovers the integer and half-integer coefficient families of Part I, is compatible with $c(0) = 1/\sqrt{\pi}$, and preserves the half-step composition requirement on the regular domain. Consequently, the complete monomial formula of the half-derivative is

$$D^{1/2}x^{\beta} = \frac{\Gamma(\beta+1)}{\Gamma\left(\beta+\frac{1}{2}\right)}x^{\beta-1/2}, \qquad x>0.$$

The Gamma function is therefore not used in this construction as an initial definition. It arises as the canonical functional form of the already-obtained algebraic coefficient law. In the next chapter, this result will be used to motivate the extension from the half-order operator to operators of general fractional order.

# Chapter II — Motivation: From Two Coefficient Chains Toward a General Order

## 2.1. The Result Obtained and the Purpose of This Chapter

In the previous chapter we continued the discrete coefficient family carried over from Part I to the corresponding real domain and obtained the monomial formula of the half-order operator

$$\boxed{D^{1/2}x^{\beta} = \frac{\Gamma(\beta+1)}{\Gamma\left(\beta+\frac{1}{2}\right)}x^{\beta-\frac{1}{2}}, \qquad x>0} \tag{2.1}$$

Formula (2.1) uniformly describes the action of the half-derivative both on integer and on half-integer powers. However, in these two cases the coefficients are expressed in different forms:

- from a monomial of integer degree the half-derivative goes to a monomial of half-integer degree, and the coefficient contains $1/\sqrt{\pi}$;
- from a monomial of half-integer degree the half-derivative goes to a monomial of integer degree, and the coefficient contains $\sqrt{\pi}$.

Thus two interlinked coefficient chains arise, whose terms alternately participate in the action of the half-order operator.

This chapter has two purposes.

The first purpose is to examine formula (2.1) in detail in two basic cases:

$$\beta = n$$

and

$$\beta = n + \frac{1}{2},$$

where $n$ is a nonnegative integer. This will clearly show how the two chains of integer and half-integer powers alternate with one another, and how the dependence on $\sqrt{\pi}$ disappears after two successive half-step actions.

The second purpose is to pose the natural question: why should we restrict ourselves only to order $\alpha = \frac{1}{2}$? The structure of the two chains suggests a natural candidate for a general $\alpha$-order operator, which will then be independently verified in the next chapter on the basis of the operator addition law.

## 2.2. The Coefficient Chain for Integer Powers

Let us consider the monomial

$$x^n, \qquad n = 0, 1, 2, \ldots.$$

Substituting $\beta = n$ in formula (2.1) we obtain

$$D^{1/2} x^n = \frac{\Gamma(n+1)}{\Gamma\left(n + \frac{1}{2}\right)} x^{n-\frac{1}{2}}$$

Since

$$\Gamma(n+1) = n!,$$

while

$$\Gamma\left(n + \frac{1}{2}\right) = \frac{(2n-1)!!}{2^n}\sqrt{\pi},$$

we have

$$\frac{\Gamma(n+1)}{\Gamma\left(n + \frac{1}{2}\right)} = \frac{n!}{\frac{(2n-1)!!}{2^n}\sqrt{\pi}} = \frac{2^n n!}{(2n-1)!!\sqrt{\pi}}$$

But

$$(2n)!! = 2^n n!,$$

therefore

$$\boxed{c(n) = \frac{(2n)!!}{(2n-1)!!\sqrt{\pi}}, \qquad n = 0, 1, 2, \ldots} \tag{2.2}$$

Accordingly,

$$\boxed{D^{1/2}x^n = \frac{(2n)!!}{(2n-1)!!\sqrt{\pi}}x^{n-\frac{1}{2}}} \tag{2.3}$$

In the case $n = 0$ the standard conventions of the double factorial are used

$$0!! = 1, \qquad (-1)!! = 1$$

therefore formula (2.2) gives

$$c\,(0) = \frac{0!!}{(-1)!!\sqrt{\pi}} = \frac{1}{\sqrt{\pi}},$$

which exactly coincides with the chosen normalization.
The first few values are:

$$D^{1/2}1 = \frac{1}{\sqrt{\pi}}x^{-1/2},$$

$$D^{1/2}x = \frac{2}{\sqrt{\pi}}x^{1/2},$$

$$D^{1/2}x^2 = \frac{8}{3\sqrt{\pi}}x^{3/2},$$

$$D^{1/2}x^3 = \frac{16}{5\sqrt{\pi}}x^{5/2},$$

$$D^{1/2}x^4 = \frac{128}{35\sqrt{\pi}}x^{7/2}$$

Thus, the chain of integer powers can be written schematically:

$$1 \mapsto \frac{1}{\sqrt{\pi}}x^{-1/2},$$

$$x \mapsto \frac{2}{\sqrt{\pi}}x^{1/2},$$

$$x^2 \mapsto \frac{8}{3\sqrt{\pi}}x^{3/2},$$

$$x^3 \mapsto \frac{16}{5\sqrt{\pi}}x^{5/2},$$

$$\dots$$

In every coefficient of this chain $\sqrt{\pi}$ occurs in the denominator.

## 2.3. The Coefficient Chain for Half-Integer Powers

Now let us consider the half-integer-power monomial

$$x^{n+\frac{1}{2}}, \qquad n = 0, 1, 2, \ldots.$$

In formula (2.1) let us write

$$\beta = n + \frac{1}{2}$$

Then

$$D^{1/2} x^{n+\frac{1}{2}} = \frac{\Gamma\left(n+\frac{3}{2}\right)}{\Gamma\left(n+1\right)} x^n \tag{2.4}$$

In the denominator we have

$$\Gamma\left(n+1\right) = n!. \tag{2.5}$$

To compute the numerator let us use the recurrence of the Gamma function:

$$\Gamma\left(n+\frac{3}{2}\right) = \left(n+\frac{1}{2}\right)\Gamma\left(n+\frac{1}{2}\right)$$

The first few values are:

$$\Gamma\left(\frac{3}{2}\right) = \frac{1}{2}\Gamma\left(\frac{1}{2}\right) = \frac{\sqrt{\pi}}{2},$$

$$\Gamma\left(\frac{5}{2}\right) = \frac{3}{2}\Gamma\left(\frac{3}{2}\right) = \frac{3\sqrt{\pi}}{4},$$

$$\Gamma\left(\frac{7}{2}\right) = \frac{5}{2}\Gamma\left(\frac{5}{2}\right) = \frac{15\sqrt{\pi}}{8},$$

$$\Gamma\left(\frac{9}{2}\right) = \frac{7}{2}\Gamma\left(\frac{7}{2}\right) = \frac{105\sqrt{\pi}}{16}$$

From this sequence we obtain the general formula

$$\boxed{\Gamma\left(n+\frac{3}{2}\right) = \frac{(2n+1)!!}{2^{n+1}}\sqrt{\pi}} \tag{2.6}$$

Formula (2.6) can also be proven by induction.
When $n = 0$,

$$\Gamma\left(\frac{3}{2}\right) = \frac{1!!}{2}\sqrt{\pi} = \frac{\sqrt{\pi}}{2}$$

Suppose that for some $n$

$$\Gamma\left(n+\frac{3}{2}\right)=\frac{(2n+1)!!}{2^{n+1}}\sqrt{\pi}$$

Then

$$\Gamma\left(n+\frac{5}{2}\right)=\left(n+\frac{3}{2}\right)\Gamma\left(n+\frac{3}{2}\right)$$

Since

$$n+\frac{3}{2}=\frac{2n+3}{2},$$

we obtain

$$\Gamma\left(n+\frac{5}{2}\right)=\frac{2n+3}{2}\cdot\frac{(2n+1)!!}{2^{n+1}}\sqrt{\pi}=\frac{(2n+3)!!}{2^{n+2}}\sqrt{\pi}$$

Thus, formula (2.6) holds for all $n\geq 0$.
Now let us substitute formulas (2.5) and (2.6) into (2.4):

$$D^{1/2}x^{n+\frac{1}{2}}=\frac{\frac{(2n+1)!!}{2^{n+1}}\sqrt{\pi}}{n!}x^n$$

Since

$$(2n)!!=2^n n!,$$

we have

$$2^{n+1}n!=2\,(2n)!!.$$

Therefore the coefficient of the half-integer power is

$$\boxed{c\left(n+\frac{1}{2}\right)=\frac{(2n+1)!!\sqrt{\pi}}{2\,(2n)!!},\qquad n=0,1,2,\dots}\tag{2.7}$$

Accordingly,

$$\boxed{D^{1/2}x^{n+\frac{1}{2}}=\frac{(2n+1)!!\sqrt{\pi}}{2\,(2n)!!}x^n}\tag{2.8}$$

The first few values are:

$$D^{1/2}x^{1/2}=\frac{\sqrt{\pi}}{2},$$

$$D^{1/2}x^{3/2}=\frac{3\sqrt{\pi}}{4}x,$$

$$D^{1/2}x^{5/2} = \frac{15\sqrt{\pi}}{16}x^2,$$

$$D^{1/2}x^{7/2} = \frac{35\sqrt{\pi}}{32}x^3$$

Thus, the chain of half-integer powers is

$$x^{1/2} \mapsto \frac{\sqrt{\pi}}{2},$$

$$x^{3/2} \mapsto \frac{3\sqrt{\pi}}{4}x,$$

$$x^{5/2} \mapsto \frac{15\sqrt{\pi}}{16}x^2,$$

$$x^{7/2} \mapsto \frac{35\sqrt{\pi}}{32}x^3,$$

$$\dots$$

In every coefficient of this chain $\sqrt{\pi}$ occurs in the numerator.

## 2.4. The Overall Picture of the Two Coefficient Chains

Let us compare the two families obtained.

For integer powers

$$\boxed{D^{1/2}x^n = \frac{(2n)!!}{(2n-1)!!\sqrt{\pi}}x^{n-\frac{1}{2}}, \qquad n = 0, 1, 2, \dots}, \tag{2.9}$$

while for half-integer powers

$$\boxed{D^{1/2}x^{n+\frac{1}{2}} = \frac{(2n+1)!!\sqrt{\pi}}{2\,(2n)!!}x^n, \qquad n = 0, 1, 2, \dots}. \tag{2.10}$$

These two formulas can be summarized by the following table.

Table 1: Integer and Half-Integer Coefficient Families

| **Original sector** | **Original monomial** | **Coefficient** | **Result of the half-derivative** | **Resulting sector** |
|---|---|---|---|---|
| Integer power | $x^n$ | $\frac{(2n)!!}{(2n-1)!!\sqrt{\pi}}$ | $\frac{(2n)!!}{(2n-1)!!\sqrt{\pi}}x^{n-\frac{1}{2}}$ | Half-integer power |
| Half-integer power | $x^{n+\frac{1}{2}}$ | $\frac{(2n+1)!!\sqrt{\pi}}{2(2n)!!}$ | $\frac{(2n+1)!!\sqrt{\pi}}{2(2n)!!}x^n$ | Integer power |

Thus, every action of $D^{1/2}$ reduces the power by exactly $\frac{1}{2}$:

$$\beta \mapsto \beta - \frac{1}{2}$$

If the original power is a positive integer, then two successive actions give the scheme

$$x^n \mapsto^{D^{1/2}} x^{n-\frac{1}{2}} \mapsto^{D^{1/2}} x^{n-1} \tag{2.11}$$

At the level of coefficients this scheme corresponds to

$$c(n) \cdot c\left(n - \frac{1}{2}\right) = n \tag{2.12}$$

The first coefficient contains $1/\sqrt{\pi}$, while the second — $\sqrt{\pi}$. Therefore in the composition these factors automatically cancel.

Table 2: The Complete Compositional Table of the First Monomials

| **Exponent** $\beta$ | **First half-step** $D^{1/2}x^{\beta}$ | **Second half-step** $D^{1/2}(D^{1/2}x^{\beta})$ | **Ordinary derivative** $Dx^{\beta}$ |
|---|---|---|---|
| $0$ | $\frac{1}{\sqrt{\pi}}x^{-1/2}$ | $0$ | $0$ |
| $\frac{1}{2}$ | $\frac{\sqrt{\pi}}{2}$ | $\frac{1}{2}x^{-1/2}$ | $\frac{1}{2}x^{-1/2}$ |
| $1$ | $\frac{2}{\sqrt{\pi}}x^{1/2}$ | $1$ | $1$ |
| $\frac{3}{2}$ | $\frac{3\sqrt{\pi}}{4}x$ | $\frac{3}{2}x^{1/2}$ | $\frac{3}{2}x^{1/2}$ |
| $2$ | $\frac{8}{3\sqrt{\pi}}x^{3/2}$ | $2x$ | $2x$ |
| $\frac{5}{2}$ | $\frac{15\sqrt{\pi}}{16}x^2$ | $\frac{5}{2}x^{3/2}$ | $\frac{5}{2}x^{3/2}$ |
| $3$ | $\frac{16}{5\sqrt{\pi}}x^{5/2}$ | $3x^2$ | $3x^2$ |

In the table the case $\beta = 0$ is special. The first action gives

$$D^{1/2}1 = \frac{1}{\sqrt{\pi}}x^{-1/2}$$

By the boundary condition obtained in Part I

$$D^{1/2}x^{-1/2} = 0,$$

therefore

$$D^{1/2}\left(D^{1/2}1\right) = 0 = D1$$

In the remaining cases the two coefficient families act alternately, and after two half-steps the ordinary first derivative is obtained exactly.

## 2.5. Compositional Verification of the Half-Integer Chain

In the previous chapter the composition condition was checked in general Gamma-function form. Now let us examine the same mechanism directly on examples of half-integer powers.

**First example:** $x^{1/2}$

By formula (2.8)

$$D^{1/2}x^{1/2} = \frac{\sqrt{\pi}}{2}$$

The result obtained is a constant. By the linearity of the operator

$$D^{1/2}\left(\frac{\sqrt{\pi}}{2}\right) = \frac{\sqrt{\pi}}{2}D^{1/2}1$$

Since

$$D^{1/2}1 = \frac{1}{\sqrt{\pi}}x^{-1/2},$$

we obtain

$$D^{1/2}\left(D^{1/2}x^{1/2}\right) = \frac{\sqrt{\pi}}{2} \cdot \frac{1}{\sqrt{\pi}}x^{-1/2}$$

Here $\sqrt{\pi}$ cancels and we are left with

$$D^{1/2}\left(D^{1/2}x^{1/2}\right) = \frac{1}{2}x^{-1/2}$$

On the other hand,

$$Dx^{1/2} = \frac{1}{2}x^{-1/2}$$

Thus,

$$\boxed{D^{1/2}\left(D^{1/2}x^{1/2}\right) = Dx^{1/2}} \tag{2.13}$$

**Second example:** $x^{3/2}$

By formula (2.8)

$$D^{1/2}x^{3/2} = \frac{3\sqrt{\pi}}{4}x$$

The second action gives

$$D^{1/2}\left(D^{1/2}x^{3/2}\right) = \frac{3\sqrt{\pi}}{4}D^{1/2}x$$

Since

$$D^{1/2}x = \frac{2}{\sqrt{\pi}}x^{1/2},$$

we obtain

$$D^{1/2}\left(D^{1/2}x^{3/2}\right) = \frac{3\sqrt{\pi}}{4}\cdot\frac{2}{\sqrt{\pi}}x^{1/2}$$

Simplifying the coefficient

$$\frac{3\sqrt{\pi}}{4}\cdot\frac{2}{\sqrt{\pi}} = \frac{3}{2}$$

Therefore

$$D^{1/2}\left(D^{1/2}x^{3/2}\right) = \frac{3}{2}x^{1/2}$$

On the other hand,

$$Dx^{3/2} = \frac{3}{2}x^{1/2}$$

Thus,

$$\boxed{D^{1/2}\left(D^{1/2}x^{3/2}\right) = Dx^{3/2}} \tag{2.14}$$

***General verification for a half-integer power***

Let us consider

$$x^{n+\frac{1}{2}}, \qquad n = 0, 1, 2, \ldots.$$

The coefficient of the first action is

$$c\left(n+\frac{1}{2}\right) = \frac{(2n+1)!!\sqrt{\pi}}{2\,(2n)!!}$$

After the first action the resulting power is $n$, whose coefficient in the second action is

$$c\,(n) = \frac{(2n)!!}{(2n-1)!!\sqrt{\pi}}$$

The product of these two coefficients is

$$c\left(n+\frac{1}{2}\right)c\,(n) = \frac{(2n+1)!!\sqrt{\pi}}{2\,(2n)!!}\cdot\frac{(2n)!!}{(2n-1)!!\sqrt{\pi}}$$

$(2n)!!$ and $\sqrt{\pi}$ cancel, and as a result we obtain

$$c\left(n+\frac{1}{2}\right)c\,(n) = \frac{(2n+1)!!}{2\,(2n-1)!!}$$

Since

$$(2n+1)!! = (2n+1)\,(2n-1)!!,$$

we have

$$c\left(n+\frac{1}{2}\right)c\,(n) = \frac{2n+1}{2} = n+\frac{1}{2} \tag{2.15}$$

Thus,

$$D^{1/2}\left(D^{1/2}x^{n+\frac{1}{2}}\right) = c\left(n+\frac{1}{2}\right)c\,(n)\,x^{n-\frac{1}{2}} \qquad = \left(n+\frac{1}{2}\right)x^{n-\frac{1}{2}} \qquad = Dx^{n+\frac{1}{2}}$$

Finally,

$$\boxed{D^{1/2}\left(D^{1/2}x^{n+\frac{1}{2}}\right) = Dx^{n+\frac{1}{2}}} \tag{2.16}$$

This verification shows that the compensation of $\sqrt{\pi}$ is not accidental. It follows directly from the basic coefficient condition

$$c\,(\beta)\,c\left(\beta-\frac{1}{2}\right) = \beta$$

## 2.6. A Natural Question: Why Only $\frac{1}{2}$ ?

The starting requirement of the half-order construction on the initial admissible monomial class was

$$D^{1/2}\circ D^{1/2} = D \tag{2.17}$$

*Here and below, formula (2.17) is understood in the restricted monomial-level sense fixed in the Introduction and in the Convention of Parameters and Notation.*

In this relation the first derivative is represented as two equal operator steps:

$$\frac{1}{2}+\frac{1}{2} = 1$$

A natural question arises: is it possible to divide the first derivative into three, four, or generally $m$ equal operator steps?

In the case of division into three equal steps we expect

$$D^{1/3}\circ D^{1/3}\circ D^{1/3} = D,$$

while for division into $m$ equal steps

$$\boxed{D^{1/m}\circ D^{1/m}\circ\cdots\circ D^{1/m} = D} \tag{2.18}$$

If we apply the operator $D^{1/m}$ successively $k$ times, we naturally obtain the operator of order $k/m$:

$$D^{1/m} \circ \cdots \circ D^{1/m} = D^{k/m} \tag{2.19}$$

Thus, at the first step positive rational orders arise

$$\alpha = \frac{k}{m}$$

The half-derivative formula

$$D^{1/2} x^{\beta} = \frac{\Gamma(\beta+1)}{\Gamma\left(\beta+1-\frac{1}{2}\right)} x^{\beta-\frac{1}{2}}$$

shows two basic changes:

- the power $\beta$ decreases by the operator order $\frac{1}{2}$;
- the argument in the denominator of the Gamma function also decreases by the same amount.

On the basis of this structure it is natural to seek, for the general $\alpha$-order operator, the form

$$D^{\alpha} x^{\beta} = C_{\alpha}(\beta)\, x^{\beta-\alpha}, \tag{2.20}$$

where the natural candidate for the coefficient is

$$\boxed{C_{\alpha}(\beta) = \frac{\Gamma(\beta+1)}{\Gamma(\beta+1-\alpha)}} \tag{2.21}$$

Accordingly, we obtain the candidate formula

$$\boxed{D^{\alpha} x^{\beta} = \frac{\Gamma(\beta+1)}{\Gamma(\beta+1-\alpha)} x^{\beta-\alpha}} \tag{2.22}$$

At this stage formula (2.22) is so far only a **motivated candidate**. The two coefficient chains and the half-order formula show why its structure is natural, but the general correctness of the formula does not yet follow from this observation alone.

In the next chapter formula (2.22) will be independently verified on the system of monomials by the requirement of the operator addition law

$$\boxed{D^{\alpha} \circ D^{\gamma} = D^{\alpha+\gamma}} \tag{2.23}$$

The formula is used in every case where the monomials and coefficients participating in the composition are defined.

First $1/3$, $1/m$ and $k/m$ order operators are considered separately. After this the general coefficient condition is obtained directly from formula (2.23), without depending on the two chains of the half-derivative.

## 2.7. Summary of the Chapter

In this chapter the Gamma-function formula of the half-derivative obtained in the previous chapter was examined on integer and half-integer powers.

For integer powers we obtained

$$D^{1/2}x^n = \frac{(2n)!!}{(2n-1)!!\sqrt{\pi}}x^{n-\frac{1}{2}},$$

while for half-integer powers

$$D^{1/2}x^{n+\frac{1}{2}} = \frac{(2n+1)!!\sqrt{\pi}}{2\,(2n)!!}x^n$$

These two coefficient families act alternately:

**Integer power → half-integer power → integer power.**

In one family $\sqrt{\pi}$ is in the denominator, while in the other it is in the numerator. Under two successive actions on admissible monomials these factors cancel, and we obtain

$$D^{1/2} \circ D^{1/2} = D$$

The half-order example naturally gave rise to the idea of dividing the first derivative into $m$ equal operator steps:

$$\left(D^{1/m}\right)^m = D$$

From this structure the candidate for the general monomial form was motivated

$$D^\alpha x^\beta = \frac{\Gamma\,(\beta+1)}{\Gamma\,(\beta+1-\alpha)}x^{\beta-\alpha}$$

In the next chapter this formula will no longer be merely a hypothesis obtained from the two chains. It will be independently verified from the operator addition law and the corresponding initial conditions.

# Chapter III — Algebraic–Operator Construction of the General Fractional Order on Monomials

## 3.1. From the Half-Derivative to the General Order

In the previous chapters, for the half-order operator we obtained the formula

$$D^{1/2}x^\beta = \frac{\Gamma\,(\beta+1)}{\Gamma\left(\beta+\frac{1}{2}\right)}x^{\beta-\frac{1}{2}}, \qquad x>0, \tag{3.1}$$

and checked the composition condition on the admissible monomial domain

$$D^{1/2} \circ D^{1/2} = D \tag{3.2}$$

Under the action of the half-derivative on a monomial the power decreases by $\frac{1}{2}$, while under two successive actions — by one:

$$\beta \mapsto \beta - \frac{1}{2} \mapsto \beta - 1$$

At the level of coefficients this process is expressed by the relation

$$c(\beta)\, c\left(\beta - \frac{1}{2}\right) = \beta. \tag{3.3}$$

Thus, $D^{1/2}$ can be regarded as each of two equal operator steps leading to the ordinary first derivative.

Now a more general question naturally arises: is it possible to divide the first derivative not only into two, but into three, four, or generally $m$ equal operator steps?

For example, we could require

$$D^{1/3} \circ D^{1/3} \circ D^{1/3} = D,$$

or, in general,

$$D^{1/m} \circ D^{1/m} \circ \cdots \circ D^{1/m} = D \tag{3.4}$$

On the basis of this idea we will first consider operators of order $1/3$ and $1/m$, then — positive rational orders $k/m$, and finally move to the general real order $\alpha$.

#### ***3.1.1. Division of the First Derivative into Three Equal Parts***

Let us seek the order-$1/3$ operator in the form

$$D^{1/3} x^{\beta} = C_{1/3}(\beta)\, x^{\beta - \frac{1}{3}}, \tag{3.5}$$

where $C_{1/3}(\beta)$ is a coefficient not yet known.

The first application of the operator gives us

$$D^{1/3} x^{\beta} = C_{1/3}(\beta)\, x^{\beta - \frac{1}{3}}$$

The second application gives

$$D^{1/3}\left(D^{1/3} x^{\beta}\right) = C_{1/3}(\beta)\, C_{1/3}\left(\beta - \frac{1}{3}\right) x^{\beta - \frac{2}{3}}$$

After the third application

$$\left(D^{1/3}\right)^3 x^{\beta} = C_{1/3}(\beta)\, C_{1/3}\left(\beta - \frac{1}{3}\right) C_{1/3}\left(\beta - \frac{2}{3}\right) x^{\beta - 1}$$

Since the threefold composition must equal the ordinary derivative,

$$\left(D^{1/3}\right)^3 x^\beta = Dx^\beta = \beta x^{\beta-1},$$

the coefficients must satisfy the condition

$$\boxed{C_{1/3}\left(\beta\right) C_{1/3}\left(\beta-\frac{1}{3}\right) C_{1/3}\left(\beta-\frac{2}{3}\right)=\beta} \tag{3.6}$$

The natural Gamma-function solution of this condition is

$$\boxed{C_{1/3}\left(\beta\right)=\frac{\Gamma\left(\beta+1\right)}{\Gamma\left(\beta+\frac{2}{3}\right)}} \tag{3.7}$$

Indeed,

$$C_{1/3}\left(\beta-\frac{1}{3}\right)=\frac{\Gamma\left(\beta+\frac{2}{3}\right)}{\Gamma\left(\beta+\frac{1}{3}\right)},$$

while

$$C_{1/3}\left(\beta-\frac{2}{3}\right)=\frac{\Gamma\left(\beta+\frac{1}{3}\right)}{\Gamma\left(\beta\right)}$$

The product of these three coefficients is

$$C_{1/3}\left(\beta\right) C_{1/3}\left(\beta-\frac{1}{3}\right) C_{1/3}\left(\beta-\frac{2}{3}\right) \qquad = \frac{\Gamma\left(\beta+1\right)}{\Gamma\left(\beta+\frac{2}{3}\right)} \cdot \frac{\Gamma\left(\beta+\frac{2}{3}\right)}{\Gamma\left(\beta+\frac{1}{3}\right)} \cdot \frac{\Gamma\left(\beta+\frac{1}{3}\right)}{\Gamma\left(\beta\right)}$$

The intermediate Gamma functions cancel and we are left with

$$\frac{\Gamma\left(\beta+1\right)}{\Gamma\left(\beta\right)}$$

By the transfer law of the Gamma function

$$\Gamma\left(\beta+1\right)=\beta\Gamma\left(\beta\right),$$

therefore

$$\frac{\Gamma\left(\beta+1\right)}{\Gamma\left(\beta\right)}=\beta$$

Thus,

$$\boxed{\left(D^{1/3}\right)^3 x^\beta = \beta x^{\beta-1} = Dx^\beta} \tag{3.8}$$

Here the same telescoping cancellation mechanism operates as in the half-derivative case; the only difference is in the number of intermediate steps.

### 3.1.2. Division of the First Derivative into $m$ Equal Parts

Now let us consider an arbitrary integer

$$m \geq 2$$

We seek an operator $D^{1/m}$ whose $m$-fold composition gives the ordinary first derivative:

$$\left(D^{1/m}\right)^m = D \tag{3.9}$$

On a monomial let us take the form

$$D^{1/m} x^\beta = C_{1/m}(\beta)\, x^{\beta - \frac{1}{m}} \tag{3.10}$$

After $m$-fold application the power changes according to the scheme

$$\beta \mapsto \beta - \frac{1}{m} \mapsto \beta - \frac{2}{m} \mapsto \cdots \mapsto \beta - \frac{m}{m} = \beta - 1$$

The product of the corresponding coefficients must satisfy the condition

$$\boxed{\prod_{j=0}^{m-1} C_{1/m}\left(\beta - \frac{j}{m}\right) = \beta} \tag{3.11}$$

Let us consider the candidate

$$\boxed{C_{1/m}(\beta) = \frac{\Gamma(\beta+1)}{\Gamma\left(\beta + 1 - \frac{1}{m}\right)}} \tag{3.12}$$

Then

$$C_{1/m}\left(\beta - \frac{j}{m}\right) = \frac{\Gamma\left(\beta + 1 - \frac{j}{m}\right)}{\Gamma\left(\beta + 1 - \frac{j+1}{m}\right)}$$

Therefore

$$\prod_{j=0}^{m-1} C_{1/m}\left(\beta - \frac{j}{m}\right) = \prod_{j=0}^{m-1} \frac{\Gamma\left(\beta + 1 - \frac{j}{m}\right)}{\Gamma\left(\beta + 1 - \frac{j+1}{m}\right)} \qquad = \frac{\Gamma(\beta+1)}{\Gamma(\beta)}$$

All intermediate Gamma functions cancel. Accordingly,

$$\prod_{j=0}^{m-1} C_{1/m}\left(\beta - \frac{j}{m}\right) = \frac{\Gamma(\beta+1)}{\Gamma(\beta)} = \beta$$

Thus,

$$\boxed{\left(D^{1/m}\right)^m x^\beta = Dx^\beta} \tag{3.13}$$

In the case $m = 2$ formula (3.12) gives

$$C_{1/2}(\beta) = \frac{\Gamma(\beta+1)}{\Gamma\left(\beta+\frac{1}{2}\right)},$$

which is exactly the coefficient of the half-derivative.

In the case $m = 3$ we obtain

$$C_{1/3}(\beta) = \frac{\Gamma(\beta+1)}{\Gamma\left(\beta+\frac{2}{3}\right)},$$

which coincides with the formula obtained in the previous subsection.

### 3.1.3. Positive Rational Order $\alpha = k/m$

If we apply $D^{1/m}$ successively $k$ times, where

$$k = 0, 1, \ldots, m,$$

the order of the resulting operator will be

$$\alpha = \frac{k}{m}$$

The composition gives

$$D^{k/m} = D^{1/m} \circ \cdots \circ D^{1/m} \tag{3.14}$$

On a monomial the telescoping product of the coefficients will be

$$C_{k/m}(\beta) = \prod_{j=0}^{k-1} C_{1/m}\left(\beta - \frac{j}{m}\right) \qquad = \prod_{j=0}^{k-1} \frac{\Gamma\left(\beta+1-\frac{j}{m}\right)}{\Gamma\left(\beta+1-\frac{j+1}{m}\right)}$$

After cancelling the intermediate terms we obtain

$$\boxed{C_{k/m}(\beta) = \frac{\Gamma(\beta+1)}{\Gamma\left(\beta+1-\frac{k}{m}\right)}} \tag{3.15}$$

Accordingly,

$$\boxed{D^{k/m}x^\beta = \frac{\Gamma(\beta+1)}{\Gamma\left(\beta+1-\frac{k}{m}\right)} x^{\beta-\frac{k}{m}}} \tag{3.16}$$

Thus, dividing the first derivative into $m$ equal steps naturally gives us the formula for positive rational orders.

The transition from rational orders to an arbitrary real $\alpha$ does not follow automatically from $m$-fold composition alone. In the following paragraphs we will regard the same form as the canonical continuation to a general real order and independently verify it on the basis of the operator addition law.

### 3.2. The Sought Form for the General-Order Operator

For the ordinary first derivative we have

$$Dx^{\beta} = \beta x^{\beta-1}$$

For the half-order we obtained

$$D^{1/2}x^{\beta} = C_{1/2}(\beta)\, x^{\beta-\frac{1}{2}}$$

For rational orders we saw in the previous paragraph that, under the action of the operator, the power of the monomial decreases exactly by the order of the operator.

Therefore for the general order $\alpha$ we naturally seek an operator of the form

$$\boxed{D^{\alpha}x^{\beta} = C_{\alpha}(\beta)\, x^{\beta-\alpha}} \tag{3.17}$$

where $C_{\alpha}(\beta)$ is a coefficient depending on $\alpha$ and $\beta$.

Formula (3.17) means that under the action of $D^{\alpha}$

$$\beta \mapsto \beta - \alpha,$$

while in front of the monomial the coefficient

$$C_{\alpha}(\beta)$$

appears.

The task is to find such a coefficient family $C_{\alpha}(\beta)$ as will:

- preserve the addition law of operator orders;

2. in the case $\alpha = 0$ give us the identity operator;
3. in the case $\alpha = 1$ give us the ordinary derivative;
4. in the case $\alpha = \frac{1}{2}$ recover the half-derivative obtained in the previous chapters.

### 3.3. The Operator Addition Law

For nonnegative integer orders repeated differentiation satisfies the law

$$D^{m} \circ D^{r} = D^{m+r}$$

For example,

$$D^2 \circ D^3 = D^5$$

For general-order operators, on the admissible monomial domain, we take as the basic requirement

$$\boxed{D^{\alpha} \circ D^{\gamma} = D^{\alpha+\gamma}} \tag{3.18}$$

Formula (3.18) means that first an action of order $\gamma$, and then an action of order $\alpha$, must equal a single action of order $\alpha + \gamma$.

This requirement is considered in the present paper as **the defining condition of the operator family built on the system of monomials**. It should not be understood as a universal law automatically fulfilled by every classical fractional operator and on every functional space.

*Formal hypothesis for (3.18): the equality is proven not simultaneously for all $\alpha, \gamma, \beta \in \mathbb{R}$, but on the domain where none of the intermediate arguments — $\beta + 1$, $\beta - \gamma + 1$, $\beta - \alpha - \gamma + 1$ — falls in the set of poles of $\Gamma$, $\{0, -1, -2, \dots\}$. This domain restriction is an inseparable hypothesis of (3.18) and of (3.22) below, and not merely an additional comment.*

The operator addition law is used only for those values of $\alpha$, $\gamma$ and $\beta$ for which all intermediate expressions participating in the composition are defined.

This reservation is not merely a formal caution; it has concrete cases in which the law is actually violated. Let us take $\beta = -\frac{1}{2}$ and $\alpha = \gamma = \frac{1}{2}$. In Section 3.16.3 it is shown that

$$D^{1/2} x^{-1/2} = 0,$$

since $\Gamma(0)$ is at a pole. Both half-steps taken separately are defined, but the result of the first step ($D^{1/2} x^{-1/2} = 0$) is an identically zero function, on which the subsequent action of $D^{1/2}$ trivially again gives zero. However, if we write the composition at the same point in a different direction ($\alpha = 1$, without $\gamma$, directly $D^1 x^{-1/2}$), we obtain

$$D^1 x^{-1/2} = -\frac{1}{2} x^{-3/2} \neq 0$$

Here $D^{1/2}\left(D^{1/2} x^{-1/2}\right)$ and $D^1 x^{-1/2}$ are unequal: the first route passes through the pole of $\Gamma(0)$ at an intermediate stage. Hence this boundary case lies outside the admissible domain of the restricted addition law and demonstrates why the identity cannot be asserted globally. The example also shows that the definedness of every intermediate expression is an essential hypothesis, rather than a merely technical qualification. The same classical counterexample is discussed by Gorenflo and Mainardi [37, Section 1.5].

## 3.4. Applying the Addition Law on a Monomial

Let us consider the monomial $x^{\beta}$.

First let us act with the order-$\gamma$ operator:

$$D^{\gamma} x^{\beta} = C_{\gamma}(\beta)\, x^{\beta-\gamma} \tag{3.19}$$

Now let us act on the obtained result with the order-$\alpha$ operator:

$$D^{\alpha}\left(D^{\gamma}x^{\beta}\right)=D^{\alpha}\left(C_{\gamma}\left(\beta\right)x^{\beta-\gamma}\right)$$

Since $C_{\gamma}(\beta)$ does not depend on $x$, by the linearity of the operator it comes out in front of the operator:

$$D^{\alpha}\left(D^{\gamma}x^{\beta}\right)=C_{\gamma}\left(\beta\right)D^{\alpha}x^{\beta-\gamma}$$

Using formula (3.17)

$$D^{\alpha}x^{\beta-\gamma}=C_{\alpha}\left(\beta-\gamma\right)x^{\beta-\gamma-\alpha}$$

Therefore

$$D^{\alpha}\left(D^{\gamma}x^{\beta}\right)=C_{\gamma}\left(\beta\right)C_{\alpha}\left(\beta-\gamma\right)x^{\beta-\alpha-\gamma} \tag{3.20}$$

On the other hand, by the addition law,

$$D^{\alpha}\left(D^{\gamma}x^{\beta}\right)=D^{\alpha+\gamma}x^{\beta}$$

In general form

$$D^{\alpha+\gamma}x^{\beta}=C_{\alpha+\gamma}\left(\beta\right)x^{\beta-\alpha-\gamma} \tag{3.21}$$

Comparing formulas (3.20) and (3.21) we obtain the basic coefficient condition

$$\boxed{C_{\gamma}\left(\beta\right)C_{\alpha}\left(\beta-\gamma\right)=C_{\alpha+\gamma}\left(\beta\right)} \tag{3.22}$$

This is the general-order coefficients' law of operator compatibility.

*This domain condition is formally an inseparable hypothesis of (3.22): the equality is proven only for those $\beta$, $\gamma$, $\alpha$ for which $\Gamma(\beta+1)$, $\Gamma(\beta-\gamma+1)$ and $\Gamma(\beta-\alpha-\gamma+1)$ are all finite and nonzero. Outside this domain (e.g., in the boundary example of §3.16.3) the composition must be checked separately by means of the reciprocal Gamma function $1/\Gamma$.*

## 3.5. The Symmetric Form of the Composition

If both orders of action are defined, then

$$D^{\alpha}\circ D^{\gamma}=D^{\alpha+\gamma}=D^{\gamma+\alpha}=D^{\gamma}\circ D^{\alpha}$$

Therefore on a monomial we have

$$C_{\gamma}\left(\beta\right)C_{\alpha}\left(\beta-\gamma\right)=C_{\alpha+\gamma}\left(\beta\right),$$

and at the same time

$$C_{\alpha}\left(\beta\right)C_{\gamma}\left(\beta-\alpha\right)=C_{\alpha+\gamma}\left(\beta\right)$$

Accordingly,

$$\boxed{C_\gamma(\beta)\, C_\alpha(\beta-\gamma) = C_\alpha(\beta)\, C_\gamma(\beta-\alpha)} \tag{3.23}$$

Formula (3.23) shows that the coefficient system must be compatible with composition of the operators in both orders.

This equality is used only in cases where all coefficients and intermediate actions on both sides are defined.

### 3.6. The Half-Derivative as a Special Case of the General Condition

In the general coefficient condition

$$C_\gamma(\beta)\, C_\alpha(\beta-\gamma) = C_{\alpha+\gamma}(\beta)$$

let us take

$$\alpha = \gamma = \frac{1}{2}$$

Then

$$\alpha + \gamma = 1,$$

and we obtain

$$C_{1/2}(\beta)\, C_{1/2}\left(\beta - \frac{1}{2}\right) = C_1(\beta) \tag{3.24}$$

Since the first-order operator is the ordinary derivative,

$$D^1 x^\beta = \beta x^{\beta-1},$$

we must have

$$C_1(\beta) = \beta$$

Thus, formula (3.24) becomes the relation

$$\boxed{C_{1/2}(\beta)\, C_{1/2}\left(\beta - \frac{1}{2}\right) = \beta} \tag{3.25}$$

This exactly coincides with the basic recurrence of the half-derivative obtained in Part I

$$c(\beta)\, c\left(\beta - \frac{1}{2}\right) = \beta$$

Accordingly,

$$C_{1/2}(\beta) = c(\beta),$$

and the general coefficient condition fully contains the half-derivative construction as a special case.

## 3.7. Requirements for Zero and Unit Orders

The general-order operator family must contain both the identity operator and the ordinary first derivative.

### *Zero order*

The zero-order operator must be the identity operator $I$:

$$D^0 = I$$

Therefore

$$D^0 x^\beta = x^\beta$$

On the other hand, in general form

$$D^0 x^\beta = C_0(\beta)\, x^\beta$$

from which it is necessary that

$$\boxed{C_0(\beta) = 1} \tag{3.26}$$

### *Unit order*

In the first-order case we must obtain the ordinary derivative:

$$D^1 x^\beta = \beta x^{\beta-1}$$

In general form

$$D^1 x^\beta = C_1(\beta)\, x^{\beta-1}$$

Therefore

$$\boxed{C_1(\beta) = \beta} \tag{3.27}$$

Thus, the initial conditions are

$$\boxed{C_0(\beta) = 1, \qquad C_1(\beta) = \beta} \tag{3.28}$$

## 3.8. The General Gamma-Function Candidate

For the half-order we obtained

$$C_{1/2}(\beta) = \frac{\Gamma(\beta+1)}{\Gamma\left(\beta+1-\frac{1}{2}\right)}$$

For rational orders as well the same structure was obtained:

$$C_{k/m}(\beta) = \frac{\Gamma(\beta+1)}{\Gamma\left(\beta+1-\frac{k}{m}\right)}$$

On the basis of these forms, for the general real order let us consider the candidate

$$\boxed{C_{\alpha}(\beta) = \frac{\Gamma(\beta+1)}{\Gamma(\beta+1-\alpha)}} \tag{3.29}$$

If this coefficient satisfies the operator condition (3.22) and the initial conditions (3.28), we obtain the action of the general-order operator on a monomial

$$\boxed{D^{\alpha}x^{\beta} = \frac{\Gamma(\beta+1)}{\Gamma(\beta+1-\alpha)}x^{\beta-\alpha}} \tag{3.30}$$

At this stage (3.29) is the canonical Gamma-function candidate motivated by rational orders. Now it must be independently verified.

## 3.9. Verification of the Gamma-Function Coefficient

Let us take

$$C_{\gamma}(\beta) = \frac{\Gamma(\beta+1)}{\Gamma(\beta+1-\gamma)} \tag{3.31}$$

For the power $\beta-\gamma$ we have

$$C_{\alpha}(\beta-\gamma) = \frac{\Gamma(\beta-\gamma+1)}{\Gamma(\beta-\gamma+1-\alpha)} \tag{3.32}$$

Now let us multiply formulas (3.31) and (3.32):

$$C_{\gamma}(\beta)\,C_{\alpha}(\beta-\gamma) = \frac{\Gamma(\beta+1)}{\Gamma(\beta+1-\gamma)} \cdot \frac{\Gamma(\beta-\gamma+1)}{\Gamma(\beta-\gamma+1-\alpha)}$$

Since

$$\beta+1-\gamma = \beta-\gamma+1,$$

the intermediate Gamma functions coincide with one another and cancel:

$$C_\gamma(\beta)\, C_\alpha(\beta-\gamma) = \frac{\Gamma(\beta+1)}{\Gamma(\beta+1-\alpha-\gamma)}$$

But by formula (3.29),

$$C_{\alpha+\gamma}(\beta) = \frac{\Gamma(\beta+1)}{\Gamma(\beta+1-(\alpha+\gamma))}$$

Thus,

$$\boxed{C_\gamma(\beta)\, C_\alpha(\beta-\gamma) = C_{\alpha+\gamma}(\beta)} \tag{3.33}$$

Accordingly,

$$\boxed{D^\alpha\left(D^\gamma x^\beta\right) = D^{\alpha+\gamma} x^\beta} \tag{3.34}$$

The formula holds in every case where both stages of the composition and the corresponding Gamma-function coefficients are defined.

By this verification it is established that the Gamma-function candidate satisfies the requirement of operator addition. Without additional conditions we do not here claim the absolute uniqueness of every possible coefficient family; formula (3.29) is the canonical solution compatible with the results of the previous chapters.

## 3.10. Verification of the Initial Conditions

Let us check the formula

$$C_\alpha(\beta) = \frac{\Gamma(\beta+1)}{\Gamma(\beta+1-\alpha)}$$

for the zero and unit orders.

***Zero order***

If

$$\alpha = 0,$$

then

$$C_0(\beta) = \frac{\Gamma(\beta+1)}{\Gamma(\beta+1)} = 1$$

Therefore

$$D^0 x^\beta = x^\beta$$

Thus,

$$\boxed{D^0 = I} \tag{3.35}$$

***Unit order***

If

$$\alpha = 1,$$

then

$$C_1(\beta) = \frac{\Gamma(\beta+1)}{\Gamma(\beta)}$$

Since

$$\Gamma(\beta+1) = \beta\Gamma(\beta),$$

we have

$$C_1(\beta) = \beta$$

Therefore

$$D^1 x^\beta = \beta x^{\beta-1} = Dx^\beta$$

Thus,

$$\boxed{D^1 = D} \tag{3.36}$$

The Gamma-function formula correctly recovers both the identity operator at $\alpha = 0$ and the ordinary derivative at $\alpha = 1$.

## 3.11. The General Formula for a Monomial

On the basis of the reasoning obtained, for a monomial we define

$$\boxed{D^\alpha x^\beta = \frac{\Gamma(\beta+1)}{\Gamma(\beta+1-\alpha)} x^{\beta-\alpha}, \qquad x > 0} \tag{3.37}$$

The formula is used for those $\alpha$ and $\beta$ for which the corresponding expression is defined. The basic real monomial domain is

$$\beta > -1, \qquad x > 0$$

In the case of a pole of the Gamma function in the denominator, the coefficient can be interpreted by means of the reciprocal Gamma function:

$$\frac{1}{\Gamma(0)} = \frac{1}{\Gamma(-1)} = \frac{1}{\Gamma(-2)} = \cdots = 0$$

Here this notation means the zeros of $1/\Gamma(z)$, and not that $\Gamma(0)$ is an ordinary finite number.

***Integer orders***

If $\alpha = m$ is a nonnegative integer, then

$$D^m x^\beta = \frac{\Gamma(\beta+1)}{\Gamma(\beta+1-m)} x^{\beta-m} \tag{3.38}$$

By successive application of the recurrence of the Gamma function

$$\frac{\Gamma(\beta+1)}{\Gamma(\beta+1-m)} = \beta(\beta-1)\cdots(\beta-m+1)$$

Therefore

$$\boxed{D^m x^\beta = \beta(\beta-1)\cdots(\beta-m+1)\, x^{\beta-m}}. \tag{3.39}$$

If $\beta = n$ is a nonnegative integer and $m \le n$, then

$$D^m x^n = \frac{n!}{(n-m)!} x^{n-m},$$

which is the ordinary formula for repeated differentiation.

If $m > n$, then

$$n + 1 - m$$

is a non-positive integer and $\Gamma(n+1-m)$ is at a pole. Interpreting by means of $1/\Gamma$, we obtain

$$D^m x^n = 0,$$

which also coincides with ordinary repeated differentiation.

***A Combined Table of Selected Orders***

### 3.11.1. The Index-Shift Law in the Γ-Normalized Basis

For the normalized monomial introduced in Section 1.11

$$e_\beta(x) = \frac{x^\beta}{\Gamma(\beta+1)}$$

the general formula (3.37) gives us

| Order $\alpha$ | Action on a monomial $D^\alpha x^\beta$ | Interpretation |
|---|---|---|
| 0 | $x^\beta$ | Identity operator |
| $\frac{1}{3}$ | $\frac{\Gamma(\beta+1)}{\Gamma\left(\beta+\frac{2}{3}\right)}x^{\beta-\frac{1}{3}}$ | One-third step of the first derivative |
| $\frac{1}{2}$ | $\frac{\Gamma(\beta+1)}{\Gamma\left(\beta+\frac{1}{2}\right)}x^{\beta-\frac{1}{2}}$ | Half-derivative |
| $\frac{2}{3}$ | $\frac{\Gamma(\beta+1)}{\Gamma\left(\beta+\frac{1}{3}\right)}x^{\beta-\frac{2}{3}}$ | Two-thirds order |
| 1 | $\beta x^{\beta-1}$ | Ordinary first derivative |
| $-\frac{1}{2}$ | $\frac{\Gamma(\beta+1)}{\Gamma\left(\beta+\frac{3}{2}\right)}x^{\beta+\frac{1}{2}}$ | Half-order Riemann–Liouville-type integral |
| $-1$ | $\frac{x^{\beta+1}}{\beta+1}$ | First-order integral taken with lower limit 0 |

$$D^\alpha e_\beta\left(x\right)=\frac{1}{\Gamma\left(\beta+1\right)}\frac{\Gamma\left(\beta+1\right)}{\Gamma\left(\beta+1-\alpha\right)}x^{\beta-\alpha}$$

After cancellation,

$$D^\alpha e_\beta\left(x\right)=\frac{x^{\beta-\alpha}}{\Gamma\left(\beta-\alpha+1\right)}=e_{\beta-\alpha}\left(x\right)$$

Thus, for all admissible $\alpha$ and $\beta$,

$$\boxed{D^\alpha e_\beta=e_{\beta-\alpha}}$$

This notation shows that in the $\Gamma$-normalized basis the fractional operator no longer requires a separate computation of the power's coefficient: it shifts the index by $\alpha$. The composition law also becomes transparent:

$$D^\gamma D^\alpha e_\beta=e_{\beta-\alpha-\gamma}=D^{\alpha+\gamma}e_\beta,$$

provided that all intermediate expressions are defined. Returning to the unnormalized monomial basis, $x^\beta=\Gamma\left(\beta+1\right)e_\beta$, again gives us (3.37). Thus, normalization does not change the existing result; it makes its structure more transparent.

### 3.11.2. Causal Normalized Power Functions and the Connection with Ortigueira

In Ortigueira's systematic approach, the normalized power function is written in causal form [33]:

$$\frac{t^\beta}{\Gamma\left(\beta+1\right)}\varepsilon\left(t\right),$$

where $\varepsilon\left(t\right)$ is the Heaviside unit step function. The corresponding power rule is

$$\boxed{D^\alpha\left[\frac{t^\beta}{\Gamma\left(\beta+1\right)}\varepsilon\left(t\right)\right]=\frac{t^{\beta-\alpha}}{\Gamma\left(\beta-\alpha+1\right)}\varepsilon\left(t\right)}$$

On $t > 0$, where $\varepsilon(t) = 1$, this result exactly coincides with the rule obtained in the previous subsection

$$D^{\alpha} e_{\beta} = e_{\beta-\alpha}$$

Thus, the present algebraic-operator construction, at the monomial level, recovers the same Γ-normalized structure as Ortigueira's causal power rule. The starting route is different: here the starting point is the composition requirement, the recurrence, and its functional continuation; while in Ortigueira's approach the causal functions are embedded in a broader systemic theory, in which distributions, convolution, and the Laplace and Fourier transforms are essential [33–36].

Therefore the coincidence of the power formulas on $t > 0$ should not be interpreted as identity of the complete operators. In the present paper agreement is proven on monomials and their finite linear combinations; the construction of full compatibility with the causal distributional space, the convolution kernel, and the transforms remains a task for further research.

### 3.12. Example: The One-Third-Order Derivative

Let us consider the monomial

$$x^2$$

By the general formula

$$D^{1/3} x^2 = \frac{\Gamma(3)}{\Gamma\left(3 - \frac{1}{3}\right)} x^{2-\frac{1}{3}}$$

Since

$$3 - \frac{1}{3} = \frac{8}{3},$$

we obtain

$$\boxed{D^{1/3} x^2 = \frac{\Gamma(3)}{\Gamma\left(\frac{8}{3}\right)} x^{5/3}} \tag{3.40}$$

Since $\Gamma(3) = 2$, we can write

$$D^{1/3} x^2 = \frac{2}{\Gamma\left(\frac{8}{3}\right)} x^{5/3}$$

Now let us apply the same operator a second time:

$$D^{1/3}\left(D^{1/3} x^2\right) = \frac{\Gamma(3)}{\Gamma\left(\frac{8}{3}\right)} D^{1/3} x^{5/3}$$

For the power $\frac{5}{3}$ we have

$$D^{1/3}x^{5/3}=\frac{\Gamma\left(\frac{8}{3}\right)}{\Gamma\left(\frac{7}{3}\right)}x^{4/3}$$

Therefore

$$D^{1/3}\left(D^{1/3}x^2\right)=\frac{\Gamma\left(3\right)}{\Gamma\left(\frac{7}{3}\right)}x^{4/3}$$

This is

$$D^{2/3}x^2$$

In the third action

$$\left(D^{1/3}\right)^3x^2=\frac{\Gamma\left(3\right)}{\Gamma\left(\frac{7}{3}\right)}D^{1/3}x^{4/3}$$

Since

$$D^{1/3}x^{4/3}=\frac{\Gamma\left(\frac{7}{3}\right)}{\Gamma\left(2\right)}x,$$

we obtain

$$\left(D^{1/3}\right)^3x^2=\frac{\Gamma\left(3\right)}{\Gamma\left(2\right)}x$$

Now

$$\Gamma\left(3\right)=2,\qquad\qquad \Gamma\left(2\right)=1,$$

therefore

$$\boxed{\left(D^{1/3}\right)^3x^2=2x=Dx^2}\tag{3.41}$$

Thus, one complete derivative is represented in the form of three equal operator steps.

**Example: A non-standard rational order and a rational power of a monomial**

The application of the general formula is not confined to integer or half-integer powers. Let us consider a case where both the order of the operator and the degree of the original monomial are rational numbers, but of a type other than half-integer.

Let us take

$$\alpha=\frac{3}{4},\ \beta=\frac{7}{8}$$

The general formula is

$$D^{\alpha}x^{\beta} = \frac{\Gamma(\beta+1)}{\Gamma(\beta+1-\alpha)}x^{\beta-\alpha} \tag{3.42}$$

Substituting the parameters we obtain

$$D^{3/4}x^{7/8} = \frac{\Gamma\left(\frac{7}{8}+1\right)}{\Gamma\left(\frac{7}{8}+1-\frac{3}{4}\right)}x^{\frac{7}{8}-\frac{3}{4}}$$

Let us compute each argument, we obtain:

$$D^{3/4}x^{7/8} = \frac{\Gamma(15/8)}{\Gamma(9/8)}x^{1/8} \tag{3.43}$$

The approximate value of the ratio of the Gamma functions is

$$D^{3/4}x^{7/8} \approx 1.012427081\, x^{1/8}$$

Therefore in numerical form

$$D^{3/4}x^{7/8} \approx 1.012427081\, x^{1/8}$$

This example clearly demonstrates the general rule:

$$\boxed{\text{original power } \beta \mapsto \text{new power } \beta-\alpha}$$

Here

$$\frac{7}{8} \mapsto \frac{1}{8}$$

In this case the values of the Gamma functions do not resolve into ordinary factorials or simple expressions containing only $\sqrt{\pi}$. Therefore the exact answer naturally remains in the form of a ratio of Gamma functions, while for practical computation its numerical approximation is used.

Let us add a second, shorter example as well:

$$\alpha = \frac{3}{4}, \ \beta = \frac{11}{9}$$

Then

$$\boxed{D^{3/4}x^{11/9} = \frac{\Gamma\left(\frac{20}{9}\right)}{\Gamma\left(\frac{53}{36}\right)}x^{17/36}}$$

In numerical form:

$$D^{3/4}x^{11/9} \approx 1.259375498\ x^{17/36}$$

## 3.13. Example: The Two-Thirds-Order Derivative

For the same monomial

$$x^2$$

let us consider $\alpha = \frac{2}{3}$.

By the general formula

$$D^{2/3}x^2 = \frac{\Gamma(3)}{\Gamma\left(3 - \frac{2}{3}\right)}x^{2-\frac{2}{3}}$$

We obtain

$$\boxed{D^{2/3}x^2 = \frac{\Gamma(3)}{\Gamma\left(\frac{7}{3}\right)}x^{4/3}} \tag{3.44}$$

Now let us act on this result with $D^{1/3}$:

$$D^{1/3}\left(D^{2/3}x^2\right) = \frac{\Gamma(3)}{\Gamma\left(\frac{7}{3}\right)}D^{1/3}x^{4/3}$$

Since

$$D^{1/3}x^{4/3} = \frac{\Gamma\left(\frac{7}{3}\right)}{\Gamma(2)}x,$$

we obtain

$$D^{1/3}\left(D^{2/3}x^2\right) = \frac{\Gamma(3)}{\Gamma(2)}x = 2x$$

Thus,

$$\boxed{D^{1/3} \circ D^{2/3}\, x^2 = D^1x^2 = 2x} \tag{3.45}$$

Similarly,

$$D^{2/3} \circ D^{1/3}\, x^2 = 2x,$$

when both compositions are defined.

This example shows that fractional operators of different orders can also be combined into one complete derivative:

$$\frac{1}{3} + \frac{2}{3} = 1$$

## 3.14. Methodological Interpretation

The general-order construction shows that the half-derivative is not an isolated or exceptional case. It represents one member of a broader operator system.

On the admissible monomial domain specified in Sections 3.3 and 3.15, the main principle is

$$D^{\alpha} \circ D^{\gamma} = D^{\alpha+\gamma}$$

Within this principle the order of the derivative behaves as an additive quantity.

For example,

$$D^{1/2} \circ D^{1/2} = D,$$

$$D^{1/3} \circ D^{1/3} \circ D^{1/3} = D,$$

$$D^{1/3} \circ D^{2/3} = D,$$

while in general

$$\left(D^{1/m}\right)^{m} = D$$

From a methodological point of view the student sees the following sequence:

**Operator composition → coefficient condition → telescoping Gamma-form → general monomial formula.**

In this way the formula

$$D^{\alpha} x^{\beta} = \frac{\Gamma\left(\beta+1\right)}{\Gamma\left(\beta+1-\alpha\right)} x^{\beta-\alpha}$$

appears not as a pre-memorized definition, but as the natural continuation of the half-step idea, rational division, and the operator addition law.

## 3.15. The Domain of Definition and the Limits of the Construction

The formula obtained in this chapter

$$D^{\alpha} x^{\beta} = \frac{\Gamma\left(\beta+1\right)}{\Gamma\left(\beta+1-\alpha\right)} x^{\beta-\alpha}$$

must be regarded as an algebraic–operator formula built on the system of monomials.

### 1. Positive $x$

For real non-integer powers we take

$$x > 0,$$

so that $x^{\beta}$ and $x^{\beta-\alpha}$ be real and single-valued.

### *2. The Degree of the Monomial*

For direct comparison with the classical left Riemann–Liouville monomial formula, we mainly consider

$$\beta > -1$$

On broader values one could consider a Gamma-function or meromorphic continuation, though this already requires additional caution.

### *3. Intermediate Actions*

The addition law

$$D^{\alpha} D^{\gamma} x^{\beta} = D^{\alpha+\gamma} x^{\beta}$$

is used only when both successive actions are defined. Formal Gamma-function cancellation does not by itself remove the need to check the intermediate domains.

### *4. Connection with the Riemann–Liouville Form*

The formula obtained coincides with the left-sided Riemann–Liouville monomial formula with lower limit 0:

$$D^{\alpha}_{RL} x^{\beta} = \frac{\Gamma(\beta+1)}{\Gamma(\beta+1-\alpha)} x^{\beta-\alpha}$$

However, coincidence with the monomial formula does not mean that the full operator theory of Riemann–Liouville has already been built in this paper. In the present construction the following are not obtained:

- the integral kernel;
- a general mechanism of non-locality;
- a theory of broad functional spaces;
- a complete system of boundary and initial conditions;
- a convolution representation.

### *5. The Difference from the Caputo Derivative*

The formula obtained gives, for a constant,

$$D^{\alpha} 1 = \frac{x^{-\alpha}}{\Gamma(1-\alpha)} \tag{3.46}$$

This result, for positive non-integer $\alpha$, is generally not zero and coincides with the Riemann–Liouville rule.

For the Caputo derivative, however, if

$$m-1<\alpha\leq m,$$

we have

$$D_C^{\alpha}x^k=0, \qquad k=0,1,\ldots,m-1$$

Thus, the Caputo and Riemann–Liouville formulas differ on the constant and low-degree terms.

The present algebraic–operator construction, with its monomial coefficients, corresponds directly to the Riemann–Liouville-type form. Agreement with the Caputo formula is obtained only for correspondingly higher-degree monomials.

### *6. The Question of Uniqueness*

Verification of the operator addition law and the initial conditions confirms that the Gamma-function coefficient is a consistent solution. Without additional conditions of regularity, continuity, and normalization, this chapter does not prove a complete classification or absolute uniqueness of every possible operator family.

### *7. Connection with Local, Integral-Free Constructions*

Other attempts are also known in the literature to construct the fractional derivative without using an integral kernel. The best known is the "conformable" derivative, proposed by Khalil and coauthors by means of a local limit [27]:

$$D_c^{\alpha}f\left(t\right)=lim_{\varepsilon\to0}\frac{f\left(t+\varepsilon t^{1-\alpha}\right)-f\left(t\right)}{\varepsilon}, \qquad 0<\alpha\leq1$$

On a monomial this operator acts by the rule $D_c^{\alpha}x^p=px^{p-\alpha}$ and, in particular, annihilates any constant: $D_c^{\alpha}1=0$ for every $\alpha$ [27,28]. The construction of the present paper differs in principle from this property: formula (3.44) shows that

$$D^{\alpha}1=\frac{x^{-\alpha}}{\Gamma\left(1-\alpha\right)}\neq0$$

for positive non-integer $\alpha$. This difference is not accidental: in the modern literature it has repeatedly been shown that the conformable derivative structurally represents an integer-order local operator in a transformed (non-integer-dimensional) space, and not, in truth, a fractional-order, non-local operator in the sense of Riemann–Liouville or Caputo [29]. The present construction, on the contrary, is specifically built so as to exactly recover, at the monomial level, the Riemann–Liouville Gamma-function coefficient (see point 4 above) — that is, it deliberately remains within the monomial "imprint" of the classical non-local definition, though the construction itself uses local, algebraic-operator means. This distinction requires emphasis, because the common feature of both approaches — avoiding the integral kernel and the limiting process at the initial stage — creates a superficial similarity, though the resulting operators belong to essentially different families. In general, the spectrum of known definitions of the fractional derivative extends far beyond the bounds of Riemann–Liouville and Caputo — the Marchaud, Grünwald–Letnikov, Riesz, and distributed-order operators are systematically discussed in [31,32]; the present paper does not fall within this broader classification, and comparison with it remains a separate task.

## 3.16. Negative Orders and the Boundary Negative Power

The general formula can also be considered for negative orders. If

$$\alpha = -\rho, \qquad \rho > 0,$$

then

$$D^{-\rho}x^{\beta} = \frac{\Gamma(\beta+1)}{\Gamma(\beta+\rho+1)}x^{\beta+\rho} \tag{3.47}$$

This formula coincides with the monomial formula of the left-sided Riemann–Liouville fractional integral of order $\rho$ with lower limit 0.

**Note on the interpretation of an intermediate order.**

$D^{-1}f$ represents the normalized integral taken with lower limit 0 and, for a positive function, has a direct geometric interpretation as area. Higher negative integer orders, for example $D^{-2}f$ and $D^{-3}f$, are the results of repeated integration. The intermediate order $D^{-5/2}f$ does not mean "half the area" between $D^{-2}f$ and $D^{-3}f$; its intermediacy concerns the order of the operator. It represents a separate value of the fractional integral, which on monomials is defined by formula (3.47) and, on the corresponding domain, satisfies the addition law of operator orders. This expresses exactly the same principle by which $D^{1/2}f$ is not "half the slope of the tangent."

### 3.16.1. Example: $\alpha = -1$

In formula (3.45) let us take

$$\rho = 1$$

Then

$$D^{-1}x^{\beta} = \frac{\Gamma(\beta+1)}{\Gamma(\beta+2)}x^{\beta+1}$$

Since

$$\Gamma(\beta+2) = (\beta+1)\,\Gamma(\beta+1),$$

we obtain

$$\boxed{D^{-1}x^{\beta} = \frac{x^{\beta+1}}{\beta+1}, \qquad \beta > -1} \tag{3.48}$$

In particular, in the case $\beta = 2$

$$D^{-1}x^{2} = \frac{\Gamma(3)}{\Gamma(4)}x^{3} = \frac{2!}{3!}x^{3} = \frac{1}{3}x^{3}$$

Thus,

$$\boxed{D^{-1}x^2 = \frac{x^3}{3}} \tag{3.49}$$

This result coincides with the integral taken with lower limit 0:

$$I_0^1 x^2 = \int_0^x t^2\,dt = \frac{x^3}{3}$$

It should not be identified with the indefinite integral

$$\int x^2\,dx = \frac{x^3}{3} + C,$$

since the operator $D^{-1}$ here defines a specific normalized integral and does not contain a free constant.

### 3.16.2. Example: $\alpha = -\frac{1}{2}$

Let us consider the monomial $x$. By the formula

$$D^{-1/2}x = \frac{\Gamma(2)}{\Gamma\left(\frac{5}{2}\right)}x^{3/2}$$

We know that

$$\Gamma(2) = 1,$$

and

$$\Gamma\left(\frac{5}{2}\right) = \frac{3}{2}\cdot\frac{1}{2}\sqrt{\pi} = \frac{3\sqrt{\pi}}{4}$$

Therefore

$$\boxed{D^{-1/2}x = \frac{4}{3\sqrt{\pi}}x^{3/2}} \tag{3.50}$$

Let us check that the half-order derivative recovers the original monomial:

$$D^{1/2}\left(D^{-1/2}x\right) = \frac{4}{3\sqrt{\pi}}D^{1/2}x^{3/2}$$

By the formula obtained in the previous chapter

$$D^{1/2}x^{3/2} = \frac{3\sqrt{\pi}}{4}x$$

Therefore

$$D^{1/2}\left(D^{-1/2}x\right) = \frac{4}{3\sqrt{\pi}} \cdot \frac{3\sqrt{\pi}}{4}x = x$$

Thus,

$$\boxed{D^{1/2} \circ D^{-1/2}\, x = x} \tag{3.51}$$

This equality, on the given monomial, expresses the mutually inverse action of positive and negative half-orders.

### 3.16.3. The Boundary Power $\beta = -\frac{1}{2}$

In Part I, for the composition of a constant function, the boundary condition was needed

$$D^{1/2}x^{-1/2} = 0$$

Now let us see how this result arises from the general Gamma-formula.

Let us take

$$\beta = -\frac{1}{2}, \qquad \alpha = \frac{1}{2}$$

Then

$$D^{1/2}x^{-1/2} = \frac{\Gamma\left(\frac{1}{2}\right)}{\Gamma(0)}x^{-1}$$

$\Gamma(0)$ is at a pole, but the reciprocal Gamma function $1/\Gamma(z)$ vanishes at $z = 0$. Therefore

$$\frac{1}{\Gamma(0)} = 0$$

in the sense of the value of the reciprocal Gamma function.

From this

$$\boxed{D^{1/2}x^{-1/2} = 0} \tag{3.52}$$

This exactly recovers the boundary coefficient obtained in Part I.

Accordingly,

$$D^{1/2}1 = \frac{1}{\sqrt{\pi}}x^{-1/2},$$

and by the second action

$$D^{1/2}\left(D^{1/2}1\right) = \frac{1}{\sqrt{\pi}}D^{1/2}x^{-1/2} = 0 = D1$$

Here it is important that the formula must not be declared automatically defined for every negative integer power. If

$$\beta = -1, -2, -3, \ldots,$$

then $\Gamma(\beta+1)$ itself is also at a pole and a separate meromorphic or distributional treatment is needed. Such an extension goes beyond the bounds of the present paper.

### 3.16.4. Connection with the Beta Function

Formula (3.45) can be independently confirmed by direct computation of the Riemann–Liouville integral. The integral taken with lower limit 0

$$D_0 x^{-\rho} x^{\beta} \equiv \frac{1}{\Gamma(\rho)} \int_0^x (x-t)^{\rho-1} t^{\beta} dt$$

represents the ordinary Riemann–Liouville-type fractional integral. Substituting

$$t = xu, \ \ u \in [0,1],$$

the integral resolves into a scaling factor and the Beta function:

$$D_0 x^{-\rho} x^{\beta} = \frac{1}{\Gamma(\rho)} x^{\beta+\rho} \cdot B(\beta+1, \rho), \tag{3.53}$$

where, by the definition of the Beta function,

$$B(\beta+1,\rho) = \frac{\Gamma(\beta+1)\Gamma(\rho)}{\Gamma(\beta+\rho+1)}$$

Substituting this expression into (3.53), the factor $1/\Gamma(\rho)$ gives us

$$D_0 x^{-\rho} x^{\beta} = \frac{\Gamma(\beta+1)}{\Gamma(\beta+\rho+1)} x^{\beta+\rho},$$

which exactly coincides with formula (3.45). Thus, for the negative order the algebraic coefficient obtained, $\Gamma(\beta+1)/\Gamma(\beta+\rho+1)$, represents an integral history (in the form of the Beta function) already, unambiguously, integrated beforehand — it does not bypass the Riemann–Liouville memory, but rather recovers it precisely at the monomial level.

## 3.17. Summary of the Chapter

In this chapter the half-order operator was generalized to the general-order operator family acting on monomials.

First the first derivative was divided into three equal steps and we obtained

$$C_{1/3}(\beta) = \frac{\Gamma(\beta+1)}{\Gamma\left(\beta+\frac{2}{3}\right)}$$

Then for $m$ equal steps we obtained

$$C_{1/m}(\beta) = \frac{\Gamma(\beta+1)}{\Gamma\left(\beta+1-\frac{1}{m}\right)},$$

while for $k$ successive steps

$$C_{k/m}(\beta) = \frac{\Gamma(\beta+1)}{\Gamma\left(\beta+1-\frac{k}{m}\right)}$$

The general-order operator was sought in the form

$$D^{\alpha}x^{\beta} = C_{\alpha}(\beta)\, x^{\beta-\alpha}$$

On the admissible monomial domain, the requirement of the operator addition law

$$D^{\alpha} \circ D^{\gamma} = D^{\alpha+\gamma}$$

gave us the coefficient condition

$$C_{\gamma}(\beta)\, C_{\alpha}(\beta-\gamma) = C_{\alpha+\gamma}(\beta)$$

It was shown that this condition, together with the initial requirements

$$C_0(\beta) = 1, \qquad C_1(\beta) = \beta$$

is satisfied by the canonical Gamma-function coefficient

$$\boxed{C_{\alpha}(\beta) = \frac{\Gamma(\beta+1)}{\Gamma(\beta+1-\alpha)}}$$

Thus, on a monomial the general formula is

$$\boxed{D^{\alpha}x^{\beta} = \frac{\Gamma(\beta+1)}{\Gamma(\beta+1-\alpha)}x^{\beta-\alpha}, \qquad x > 0}$$

On nonnegative integer orders this formula reverts to ordinary repeated differentiation. On positive non-integer orders it coincides with the left-sided Riemann–Liouville monomial formula with lower limit 0, while on negative orders — with the monomial form of the corresponding Riemann–Liouville-type integral.

At the same time the limits of the construction were determined: the result obtained is an algebraic–operator formula built on the system of monomials, and coincidence with the classical form at the monomial level does not yet represent the construction of a complete integral, non-local, or functional-analytic theory.

In the next chapter this operator will be linearly extended to polynomial initial functions, and the interrelation of the constant term, low degrees, and the Riemann–Liouville and Caputo formulas will be considered separately.

# Chapter IV — Linear Extension of the Fractional Operator to Polynomial Initial Functions

## 4.1. From a Monomial to a Polynomial

In the previous chapter the general-order operator on a monomial was defined by the formula

$$\boxed{D^{\alpha}x^{\beta} = \frac{\Gamma\left(\beta+1\right)}{\Gamma\left(\beta+1-\alpha\right)}x^{\beta-\alpha}, \qquad x>0} \tag{4.1}$$

For a polynomial the exponents of the powers are nonnegative integers. Therefore substituting $\beta = j$ we obtain

$$\boxed{D^{\alpha}x^{j} = \frac{\Gamma\left(j+1\right)}{\Gamma\left(j+1-\alpha\right)}x^{j-\alpha}, \qquad j=0,1,2,\ldots} \tag{4.2}$$

Now a natural question arises: how should we act not only on a single monomial, but on a polynomial?

Let us take a degree-$N$ polynomial

$$P\left(x\right) = a_0 + a_1x + a_2x^2 + \cdots + a_Nx^N, \tag{4.3}$$

or, in abbreviated notation,

$$P\left(x\right) = \sum_{j=0}^{N} a_jx^j \tag{4.4}$$

The ordinary derivative is a linear operator:

$$D\left(f+g\right) = Df + Dg, \tag{4.5}$$

and

$$D\left(\lambda f\right) = \lambda Df, \tag{4.6}$$

where $\lambda$ is a constant number.

Therefore the ordinary derivative acts term-by-term on a polynomial:

$$DP\left(x\right) = \sum_{j=0}^{N} a_jD\left(x^j\right) = \sum_{j=1}^{N} j\,a_jx^{j-1} \tag{4.7}$$

By the same principle we will also regard the general-order operator as a linear operator, and extend the action already defined on monomials term-by-term to polynomial initial functions.

Here the words "extension to a polynomial" mean that the initial function of the operator is a polynomial. In the case of non-integer $\alpha$ the result obtained will, in general, no longer be an ordinary polynomial, since its terms will have degrees

$$j-\alpha$$

This issue is discussed in detail in §4.4.

## 4.2. The Principle of Linearity

For the general-order operator let us require linearity:

$$\boxed{D^{\alpha}(f+g)=D^{\alpha}f+D^{\alpha}g} \tag{4.8}$$

and

$$\boxed{D^{\alpha}(\lambda f)=\lambda D^{\alpha}f} \tag{4.9}$$

Here $\lambda$ is a constant coefficient independent of $x$.

The requirement of linearity is natural for two reasons.

First, the ordinary derivative satisfies exactly this property. If $D^{\alpha}$ is regarded as a continuation of the differentiation operator family, it is desirable that it preserve this basic algebraic property.

Second, the action defined on monomials can be extended to their finite linear combinations precisely by means of linearity.

If

$$P(x)=\sum_{j=0}^{N}a_jx^j,$$

then we define

$$D^{\alpha}P(x)=D^{\alpha}\left(\sum_{j=0}^{N}a_jx^j\right)$$

On the basis of linearity

$$D^{\alpha}P(x)=\sum_{j=0}^{N}D^{\alpha}\left(a_jx^j\right),$$

while, taking the constant coefficients out,

$$\boxed{D^{\alpha}P(x)=\sum_{j=0}^{N}a_jD^{\alpha}x^j} \tag{4.10}$$

Thus, to determine the action on a polynomial it is sufficient to apply formula (4.2) to each monomial separately.

## 4.3. The General Formula for a Polynomial Initial Function

For a monomial we have

$$D^{\alpha}x^{j}=\frac{\Gamma\left(j+1\right)}{\Gamma\left(j+1-\alpha\right)}x^{j-\alpha}$$

Therefore formula (4.10) gives us

$$D^{\alpha}P\left(x\right)=\sum_{j=0}^{N}a_{j}\frac{\Gamma\left(j+1\right)}{\Gamma\left(j+1-\alpha\right)}x^{j-\alpha}$$

Thus, if

$$P\left(x\right)=\sum_{j=0}^{N}a_{j}x^{j},$$

then

$$\boxed{D^{\alpha}P\left(x\right)=\sum_{j=0}^{N}a_{j}\frac{\Gamma\left(j+1\right)}{\Gamma\left(j+1-\alpha\right)}x^{j-\alpha},\qquad x>0}\tag{4.11}$$

This is the linear extension, to polynomial initial functions, of the operator built on monomials. Since

$$x^{j-\alpha}=x^{-\alpha}x^{j},$$

formula (4.11) can also be written in another form:

$$\boxed{D^{\alpha}P\left(x\right)=x^{-\alpha}\sum_{j=0}^{N}a_{j}\frac{\Gamma\left(j+1\right)}{\Gamma\left(j+1-\alpha\right)}x^{j}}\tag{4.12}$$

This writing shows us that all terms have the common factor $x^{-\alpha}$, while in the brackets remains a finite sum of integer powers of $x$ with the corresponding Gamma-function coefficients.

If $\alpha$ is a nonnegative integer, formula (4.11) reverts to the ordinary rule of repeated differentiation. If $\alpha$ is non-integer, the resulting function is a finite sum of generalized power terms.

## 4.4. The Functional Form of the Result

In the case of non-integer $\alpha$

$$D^{\alpha}P\left(x\right)$$

in general no longer belongs to the ordinary space of polynomials $\mathbb{R}[x]$.

Indeed, the degrees of the terms of the original polynomial are

$$0, 1, 2, \ldots, N,$$

while after the action of the operator they turn into degrees

$$-\alpha, \qquad 1-\alpha, \qquad 2-\alpha, \qquad \ldots, \qquad N-\alpha$$

Therefore

$$D^{\alpha} : span\{1, x, \ldots, x^N\} \to span\{x^{-\alpha}, x^{1-\alpha}, \ldots, x^{N-\alpha}\}. \tag{4.13}$$

This means that the operator transforms the polynomial initial space into a shifted power-space.

For example, for the half-order

$$D^{1/2} : span\{1, x, x^2, \ldots, x^N\} \to span\left\{x^{-1/2}, x^{1/2}, x^{3/2}, \ldots, x^{N-1/2}\right\} \tag{4.14}$$

Therefore a more precise formulation is not "the fractional derivative is again a polynomial," but rather:

the fractional operator defined on monomials acts linearly on polynomial initial functions and transforms them into finite sums of shifted powers.

If the polynomial has a nonzero constant term, in the resulting expression the term $x^{-\alpha}$ appears, which can be singular at $x = 0$. It is precisely for this reason that the whole chapter is considered on the domain

$$x > 0$$

### 4.5. An Important Remark on the Constant Term

For the ordinary derivative the derivative of the constant term is zero:

$$D(a_0) = 0 \tag{4.15}$$

For the general-order monomial operator, however, the constant must be regarded as

$$a_0 = a_0 x^0$$

On the basis of linearity

$$D^{\alpha} a_0 = a_0 D^{\alpha} x^0 \tag{4.16}$$

Substituting $j = 0$ in formula (4.2) we obtain

$$D^{\alpha} x^0 = \frac{\Gamma(1)}{\Gamma(1-\alpha)} x^{-\alpha}$$

Since

$$\Gamma(1) = 1,$$

we have

$$\boxed{D^{\alpha} 1 = \frac{x^{-\alpha}}{\Gamma(1-\alpha)}} \tag{4.17}$$

Therefore

$$\boxed{D^{\alpha} a_0 = \frac{a_0}{\Gamma(1-\alpha)} x^{-\alpha}} \tag{4.18}$$

For positive non-integer $\alpha$ this expression is generally not zero.

In the case of the half-derivative

$$\alpha = \frac{1}{2},$$

therefore

$$D^{1/2} a_0 = \frac{a_0}{\Gamma\left(\frac{1}{2}\right)} x^{-1/2}$$

Since

$$\Gamma\left(\frac{1}{2}\right) = \sqrt{\pi},$$

we obtain

$$\boxed{D^{1/2} a_0 = \frac{a_0}{\sqrt{\pi}} x^{-1/2}} \tag{4.19}$$

This is an important difference between the ordinary and the non-integer-order Riemann–Liouville-type monomial operators.

If $\alpha = m$ is a positive integer, then

$$1 - \alpha = 1 - m$$

is a non-positive integer and $\Gamma(1-m)$ is at a pole. By the value of the reciprocal Gamma function

$$\frac{1}{\Gamma(1-m)} = 0,$$

therefore

$$D^{m} a_0 = 0,$$

which coincides with ordinary repeated differentiation.

### 4.6. The Difference Between Riemann–Liouville and Caputo on Polynomials

Formula (4.11) exactly coincides with the monomial rule of the left-sided Riemann–Liouville derivative with lower limit 0. Therefore for a polynomial the Riemann–Liouville form is

$$\boxed{D^{\alpha}_{RL} P(x) = \sum_{j=0}^{N} a_j \frac{\Gamma(j+1)}{\Gamma(j+1-\alpha)} x^{j-\alpha}} \tag{4.20}$$

The Caputo derivative [20], which is often used in initial-value problems, acts differently on low-degree terms (see also [21–24] for a comparative review).

Suppose

$$m-1 < \alpha < m, \qquad m \in \mathbb{N}.$$

Then the Caputo derivative annihilates every term whose degree is less than $m$:

$$D^{\alpha}_{C} x^{j} = 0, \qquad j = 0, 1, \ldots, m-1 \tag{4.21}$$

while for $j \geq m$

$$D^{\alpha}_{C} x^{j} = \frac{\Gamma(j+1)}{\Gamma(j+1-\alpha)} x^{j-\alpha} \tag{4.22}$$

Therefore for a polynomial

$$\boxed{D^{\alpha}_{C} P(x) = \sum_{j=m}^{N} a_j \frac{\Gamma(j+1)}{\Gamma(j+1-\alpha)} x^{j-\alpha}} \tag{4.23}$$

If $N < m$, this sum equals zero.

The difference between the Riemann–Liouville and Caputo results is

$$\boxed{D^{\alpha}_{RL} P(x) - D^{\alpha}_{C} P(x) = \sum_{j=0}^{m-1} a_j \frac{\Gamma(j+1)}{\Gamma(j+1-\alpha)} x^{j-\alpha}} \tag{4.24}$$

Thus, the difference is created only by the low-degree part of the polynomial.

In particular, if

$$0 < \alpha < 1,$$

then $m = 1$, therefore the Caputo derivative annihilates only the constant term:

$$D^{\alpha}_{C} a_0 = 0$$

In this case

$$\boxed{D_{RL}^{\alpha}P\left(x\right)-D_{C}^{\alpha}P\left(x\right)=\frac{a_0}{\Gamma\left(1-\alpha\right)}x^{-\alpha}} \tag{4.25}$$

If

$$P\left(0\right)=a_0=0,$$

then for $0<\alpha<1$ the two formulas coincide on the polynomial:

$$D_{RL}^{\alpha}P\left(x\right)=D_{C}^{\alpha}P\left(x\right) \tag{4.26}$$

Table 3: The action of the Riemann–Liouville and Caputo operators for $0<\alpha<1$

| **Initial function** | **Riemann–Liouville result** | **Caputo result** |
|---|---|---|
| 1 | $\frac{x^{-\alpha}}{\Gamma(1-\alpha)}$ | 0 |
| $a_0$ | $\frac{a_0 x^{-\alpha}}{\Gamma(1-\alpha)}$ | 0 |
| $x^j$, $j\geq 1$ | $\frac{\Gamma(j+1)}{\Gamma(j+1-\alpha)}x^{j-\alpha}$ | same |
| $P\left(x\right)$, $P\left(0\right)=0$ | Term-by-term formula obtained | same |
| General $P\left(x\right)$ | Includes the contribution of the constant term | Constant term vanishes |

The present algebraic–operator construction follows the Riemann–Liouville rule with respect to the constant term. Agreement with the Caputo form is obtained only on the corresponding terms.

## 4.7. Example I — A Linear Polynomial

Let us take

$$P\left(x\right)=a+bx \tag{4.27}$$

Its half-derivative, by linearity, is

$$D^{1/2}P\left(x\right)=aD^{1/2}1+bD^{1/2}x \tag{4.28}$$

We already know that

$$D^{1/2}1=\frac{1}{\sqrt{\pi}}x^{-1/2}, \tag{4.29}$$

and

$$D^{1/2}x=\frac{2}{\sqrt{\pi}}x^{1/2} \tag{4.30}$$

Therefore

$$D^{1/2}\left(a+bx\right)=\frac{a}{\sqrt{\pi}}x^{-1/2}+\frac{2b}{\sqrt{\pi}}x^{1/2}$$

Thus,

$$\boxed{D^{1/2}\left(a+bx\right)=\frac{1}{\sqrt{\pi}}\left(ax^{-1/2}+2bx^{1/2}\right)} \tag{4.31}$$

This is the Riemann–Liouville-type half-derivative of a linear polynomial.

For the Caputo half-derivative the constant term vanishes and we obtain

$$\boxed{D_C^{1/2}\left(a+bx\right)=\frac{2b}{\sqrt{\pi}}x^{1/2}} \tag{4.32}$$

The difference between the two results is

$$\frac{a}{\sqrt{\pi}}x^{-1/2}$$

## 4.8. Example II — A Quadratic Polynomial

Let us consider

$$P\left(x\right)=a+bx+cx^2 \tag{4.33}$$

By linearity

$$D^{1/2}P\left(x\right)=aD^{1/2}1+bD^{1/2}x+cD^{1/2}x^2 \tag{4.34}$$

We know:

$$D^{1/2}1=\frac{1}{\sqrt{\pi}}x^{-1/2},$$

$$D^{1/2}x=\frac{2}{\sqrt{\pi}}x^{1/2},$$

and

$$D^{1/2}x^2=\frac{8}{3\sqrt{\pi}}x^{3/2}$$

Therefore

$$D^{1/2}\left(a+bx+cx^2\right)=\frac{a}{\sqrt{\pi}}x^{-1/2}+\frac{2b}{\sqrt{\pi}}x^{1/2}+\frac{8c}{3\sqrt{\pi}}x^{3/2}$$

Taking out the common factor

$$\boxed{D^{1/2}\left(a+bx+cx^2\right)=\frac{1}{\sqrt{\pi}}\left(ax^{-1/2}+2bx^{1/2}+\frac{8c}{3}x^{3/2}\right)} \tag{4.35}$$

For the Caputo half-derivative we obtain

$$\boxed{D_C^{1/2}\left(a + bx + cx^2\right) = \frac{1}{\sqrt{\pi}}\left(2bx^{1/2} + \frac{8c}{3}x^{3/2}\right)} \tag{4.36}$$

## 4.9. Example III — A Concrete Polynomial

Let us take the concrete polynomial

$$P(x) = 1 + 3x + 2x^2 \tag{4.37}$$

Let us find its half-derivative.

By linearity

$$D^{1/2}P(x) = D^{1/2}1 + 3D^{1/2}x + 2D^{1/2}x^2 \tag{4.38}$$

Let us use the known formulas:

$$D^{1/2}1 = \frac{1}{\sqrt{\pi}}x^{-1/2},$$

$$D^{1/2}x = \frac{2}{\sqrt{\pi}}x^{1/2},$$

$$D^{1/2}x^2 = \frac{8}{3\sqrt{\pi}}x^{3/2}$$

Then

$$D^{1/2}P(x) = \frac{1}{\sqrt{\pi}}x^{-1/2} + 3\cdot\frac{2}{\sqrt{\pi}}x^{1/2} + 2\cdot\frac{8}{3\sqrt{\pi}}x^{3/2} \qquad = \frac{1}{\sqrt{\pi}}x^{-1/2} + \frac{6}{\sqrt{\pi}}x^{1/2} + \frac{16}{3\sqrt{\pi}}x^{3/2}$$

Thus,

$$\boxed{D^{1/2}\left(1 + 3x + 2x^2\right) = \frac{1}{\sqrt{\pi}}\left(x^{-1/2} + 6x^{1/2} + \frac{16}{3}x^{3/2}\right)} \tag{4.39}$$

The Caputo half-derivative will be

$$\boxed{D_C^{1/2}\left(1 + 3x + 2x^2\right) = \frac{1}{\sqrt{\pi}}\left(6x^{1/2} + \frac{16}{3}x^{3/2}\right)} \tag{4.40}$$

## 4.10. Verification of the Double Half-Step Action

Now let us check the main compositional condition on the initial polynomial class:

$$D^{1/2}\left(D^{1/2}P\right) = DP \tag{4.41}$$

Let us again take

$$P(x) = 1 + 3x + 2x^2$$

Its ordinary derivative is

$$DP(x) = 3 + 4x \tag{4.42}$$

By the first half-step action we obtained

$$D^{1/2}P(x) = \frac{1}{\sqrt{\pi}}x^{-1/2} + \frac{6}{\sqrt{\pi}}x^{1/2} + \frac{16}{3\sqrt{\pi}}x^{3/2} \tag{4.43}$$

Now let us apply the operator a second time, term-by-term.

***The term arising from the constant term***

For the first term we need

$$D^{1/2}x^{-1/2}$$

By the general formula

$$D^{1/2}x^{-1/2} = \frac{\Gamma\left(\frac{1}{2}\right)}{\Gamma(0)}x^{-1}$$

$\Gamma(0)$ is at a pole, while the reciprocal Gamma function $1/\Gamma(z)$ vanishes at $z = 0$. Therefore

$$\boxed{D^{1/2}x^{-1/2} = 0} \tag{4.44}$$

Accordingly,

$$D^{1/2}\left(\frac{1}{\sqrt{\pi}}x^{-1/2}\right) = 0 \tag{4.45}$$

This coincides with the ordinary derivative of the constant term.

***The term arising from the linear term***

We know that

$$D^{1/2}x^{1/2} = \frac{\sqrt{\pi}}{2}$$

Therefore

$$D^{1/2}\left(\frac{6}{\sqrt{\pi}}x^{1/2}\right) = \frac{6}{\sqrt{\pi}} \cdot \frac{\sqrt{\pi}}{2} \qquad = 3 \tag{4.46}$$

This coincides with

$$D\left(3x\right)=3$$

***The term arising from the quadratic term***

We know that

$$D^{1/2}x^{3/2}=\frac{3\sqrt{\pi}}{4}x$$

Therefore

$$D^{1/2}\left(\frac{16}{3\sqrt{\pi}}x^{3/2}\right)=\frac{16}{3\sqrt{\pi}}\cdot\frac{3\sqrt{\pi}}{4}x \qquad =4x \tag{4.47}$$

This coincides with

$$D\left(2x^{2}\right)=4x$$

Adding all three results

$$D^{1/2}\left(D^{1/2}P\left(x\right)\right)=0+3+4x$$

Thus,

$$\boxed{D^{1/2}\left(D^{1/2}\left(1+3x+2x^{2}\right)\right)=3+4x=D\left(1+3x+2x^{2}\right)} \tag{4.48}$$

## 4.11. A General Verification for a Polynomial

Let us take the general polynomial

$$P\left(x\right)=\sum_{j=0}^{N}a_{j}x^{j} \tag{4.49}$$

Its order-$\alpha$ action is

$$D^{\alpha}P\left(x\right)=\sum_{j=0}^{N}a_{j}\frac{\Gamma\left(j+1\right)}{\Gamma\left(j+1-\alpha\right)}x^{j-\alpha} \tag{4.50}$$

Now let us act on the obtained result with the order-$\gamma$ operator:

$$D^{\gamma}\left(D^{\alpha}P\left(x\right)\right)=\sum_{j=0}^{N}a_{j}\frac{\Gamma\left(j+1\right)}{\Gamma\left(j+1-\alpha\right)}D^{\gamma}x^{j-\alpha} \tag{4.51}$$

By the monomial formula

$$D^{\gamma}x^{j-\alpha}=\frac{\Gamma\left(j-\alpha+1\right)}{\Gamma\left(j-\alpha+1-\gamma\right)}x^{j-\alpha-\gamma} \tag{4.52}$$

Therefore

$$D^{\gamma}\left(D^{\alpha}P\left(x\right)\right)=\sum_{j=0}^{N}a_{j}\frac{\Gamma\left(j+1\right)}{\Gamma\left(j+1-\alpha\right)}\frac{\Gamma\left(j+1-\alpha\right)}{\Gamma\left(j+1-\alpha-\gamma\right)}x^{j-\alpha-\gamma}$$

Cancelling the intermediate Gamma functions we obtain

$$D^{\gamma}\left(D^{\alpha}P\left(x\right)\right)=\sum_{j=0}^{N}a_{j}\frac{\Gamma\left(j+1\right)}{\Gamma\left(j+1-\alpha-\gamma\right)}x^{j-\alpha-\gamma} \tag{4.53}$$

On the other hand,

$$D^{\alpha+\gamma}P\left(x\right)=\sum_{j=0}^{N}a_{j}\frac{\Gamma\left(j+1\right)}{\Gamma\left(j+1-\left(\alpha+\gamma\right)\right)}x^{j-\alpha-\gamma} \tag{4.54}$$

The denominators are identical:

$$j+1-\alpha-\gamma=j+1-\left(\alpha+\gamma\right)$$

Thus,

$$\boxed{D^{\gamma}\left(D^{\alpha}P\left(x\right)\right)=D^{\alpha+\gamma}P\left(x\right)} \tag{4.55}$$

Formula (4.55) holds only when every intermediate monomial action occurring term-by-term is defined on the admissible domain.

Formula (4.55) does not represent an independent new proof for polynomials. It is the direct result of the addition law and linearity already obtained on monomials.

In particular, if

$$\alpha=\gamma=\frac{1}{2},$$

we obtain

$$\boxed{D^{1/2}\left(D^{1/2}P\right)=DP} \tag{4.56}$$

*Thus, formula (4.56) is a statement about the initial polynomial class (with the constant-term boundary action interpreted by the reciprocal Gamma function); it does not imply the same identity on the entire enlarged graded space.*

For the constant term this equality rests on the boundary value

$$D^{1/2}x^{-1/2}=0$$

## 4.12. A Summary Table of Formulas

The table shows that the half-derivative transforms each term independently, while after the second half-step the ordinary derivative is obtained.

Table 4: The Half-Derivative on Polynomial Initial Functions

| **Initial function** $f(x)$ | $D^{1/2}f(x)$ | **Second half-step** | $Df(x)$ |
|---|---|---|---|
| $a$ | $\frac{a}{\sqrt{\pi}}x^{-1/2}$ | $0$ | $0$ |
| $bx$ | $\frac{2b}{\sqrt{\pi}}x^{1/2}$ | $b$ | $b$ |
| $cx^2$ | $\frac{8c}{3\sqrt{\pi}}x^{3/2}$ | $2cx$ | $2cx$ |
| $dx^3$ | $\frac{16d}{5\sqrt{\pi}}x^{5/2}$ | $3dx^2$ | $3dx^2$ |
| $a+bx$ | $\frac{a}{\sqrt{\pi}}x^{-1/2}+\frac{2b}{\sqrt{\pi}}x^{1/2}$ | $b$ | $b$ |
| $a+bx+cx^2$ | $\frac{a}{\sqrt{\pi}}x^{-1/2}+\frac{2b}{\sqrt{\pi}}x^{1/2}+\frac{8c}{3\sqrt{\pi}}x^{3/2}$ | $b+2cx$ | $b+2cx$ |

## 4.13. Methodological Significance

Polynomials are particularly convenient for introducing the idea of the general-order operator, for several reasons.

First, a polynomial is a finite sum of monomials. Therefore it is sufficient for the student to understand how the operator acts on one monomial; the result then naturally extends to the whole polynomial by means of linearity.

Second, the role of the coefficients is clearly visible on polynomials. Each term has its own degree $j$, which determines the corresponding Gamma-function coefficient

$$\frac{\Gamma(j+1)}{\Gamma(j+1-\alpha)}$$

appears in front of it.

Third, the difference between the ordinary and the non-integer-order operators is clearly visible. The ordinary derivative annihilates the constant term and transforms a polynomial again into a polynomial. The non-integer-order Riemann–Liouville-type action, however, produces $x^{-\alpha}$ from the constant term and, in general, moves the result into the space of shifted powers.

Fourth, on polynomials it is easy to compare the Riemann–Liouville and Caputo results. The difference is entirely concentrated on the low-degree terms.

Fifth, polynomials are algebraic objects familiar to the student. Therefore the idea of the fractional order does not enter from the outset through complex integral kernels or abstract functional spaces. It begins with term-by-term computations, and only afterward is connected to the classical theory.

Thus, polynomials create a natural learning environment where the idea of the general-order operator can be unfolded step by step.

## 4.14. Connection with the Previous Construction

The approach used in this chapter fully continues the logic of the previous chapters.

Originally we had the ordinary derivative on a monomial:

$$Dx^n = nx^{n-1}$$

Then we considered the half-derivative:

$$D^{1/2}x^{\beta} = c(\beta)\, x^{\beta-\frac{1}{2}},$$

and for the coefficients we obtained the condition

$$c(\beta)\, c\left(\beta - \frac{1}{2}\right) = \beta$$

After functional continuation we obtained

$$c(\beta) = \frac{\Gamma(\beta+1)}{\Gamma\left(\beta+\frac{1}{2}\right)}$$

Then the construction was generalized for any corresponding order:

$$D^{\alpha}x^{\beta} = \frac{\Gamma(\beta+1)}{\Gamma(\beta+1-\alpha)} x^{\beta-\alpha}$$

In this chapter the monomial formula was linearly extended to polynomial initial functions:

$$\boxed{D^{\alpha}\left(\sum_{j=0}^{N} a_j x^j\right) = \sum_{j=0}^{N} a_j \frac{\Gamma(j+1)}{\Gamma(j+1-\alpha)} x^{j-\alpha}} \tag{4.57}$$

With this, the algebraic–operator construction moved from a single monomial to a finite-dimensional class of functions.

## 4.15. Summary of the Chapter and Transition to the Next Chapter

In this chapter the general-order operator defined on monomials was linearly extended to polynomial initial functions.

If

$$P(x) = \sum_{j=0}^{N} a_j x^j,$$

then we obtained

$$\boxed{D^{\alpha}P(x) = \sum_{j=0}^{N} a_j \frac{\Gamma(j+1)}{\Gamma(j+1-\alpha)} x^{j-\alpha}}$$

It was emphasized that in the case of non-integer $\alpha$ the result is, in general, no longer a polynomial. It is a finite sum of shifted powers.

The constant term was considered separately:

$$D^{\alpha}a_0 = \frac{a_0}{\Gamma(1-\alpha)} x^{-\alpha}$$

This formula corresponds to the Riemann–Liouville monomial rule. The Caputo operator annihilates the corresponding low-degree terms.

For the half-order, on the examples and on the initial polynomial class, the compositional condition was checked

$$D^{1/2}\left(D^{1/2}P\right) = DP$$

More generally, on the corresponding admissible term-by-term domain, we obtained

$$D^{\gamma}\left(D^{\alpha}P\right) = D^{\alpha+\gamma}P$$

In the next chapter the formulas obtained will be presented with computational and visual examples. We will examine the initial graphs of monomials and polynomials, their half-derivatives, and their ordinary derivatives, in order to compare the operator steps visually with one another.

# Chapter V — Computational and Visual Examples

## 5.1. The Purpose of the Chapter

In the previous chapters, for the general-order operator on a monomial we obtained the formula

$$\boxed{D^{\alpha}x^{\beta} = \frac{\Gamma\left(\beta+1\right)}{\Gamma\left(\beta+1-\alpha\right)}x^{\beta-\alpha}, \qquad x>0} \tag{5.1}$$

while for the half-order case

$$\boxed{D^{1/2}x^{\beta} = \frac{\Gamma\left(\beta+1\right)}{\Gamma\left(\beta+\frac{1}{2}\right)}x^{\beta-\frac{1}{2}}} \tag{5.2}$$

We also showed that, on the corresponding domain of definition,

$$\boxed{D^{1/2}\left(D^{1/2}x^{\beta}\right) = Dx^{\beta}} \tag{5.3}$$

For polynomial initial functions the operator was defined term-by-term. If

$$P\left(x\right) = \sum_{j=0}^{N} a_j x^j,$$

then

$$\boxed{D^{\alpha}P\left(x\right) = \sum_{j=0}^{N} a_j \frac{\Gamma\left(j+1\right)}{\Gamma\left(j+1-\alpha\right)}x^{j-\alpha}} \tag{5.4}$$

The purpose of the present chapter is not to re-prove these formulas. Their algebraic justification has already been given in the previous chapters. Here we will examine, for the results obtained:

- concrete computational examples;
- numerical verification;
- tabular representation;
- graphical interpretation;
- the possibility of computer visualization.

In the half-order case we will compare three functions:

$$f(x), \qquad D^{1/2} f(x), \qquad Df(x) \tag{5.5}$$

For the admissible monomial and polynomial examples considered below, we will additionally check that

$$D^{1/2}\left(D^{1/2} f\right) = Df \tag{5.6}$$

The term "intermediate" is used here in the sense of **the order of the operator**: $D^{1/2}$ is a half-order action lying between the zero- and first-order operators. This does not mean that the numerical value of $D^{1/2} f(x)$ or its graph necessarily lies, at every point, between $f(x)$ and $Df(x)$.

## 5.2. The General Scheme for Visual Comparison

For the monomial

$$f_\beta(x) = x^\beta$$

let us define its half-derivative

$$g_\beta(x) = D^{1/2} f_\beta(x) = \frac{\Gamma(\beta+1)}{\Gamma\left(\beta+\frac{1}{2}\right)} x^{\beta-\frac{1}{2}} \tag{5.7}$$

The second half-step action will be

$$h_\beta(x) = D^{1/2} g_\beta(x) \tag{5.8}$$

By the compositional law,

$$h_\beta(x) = \beta x^{\beta-1} = Df_\beta(x) \tag{5.9}$$

Thus, graphically we can compare:

$$f_\beta(x) = x^\beta,$$

$$g_\beta(x) = \frac{\Gamma(\beta+1)}{\Gamma\left(\beta+\frac{1}{2}\right)} x^{\beta-\frac{1}{2}},$$

and

$$h_\beta(x) = Df_\beta(x) = \beta x^{\beta-1}$$

Because of the fractional powers and possible singularity, the graphs must be built on a positive interval

$$x \in [\varepsilon, L], \qquad \varepsilon > 0 \tag{5.10}$$

For example, we could take

$$\varepsilon = 0.01, \qquad L = 5$$

## 5.3. Example I — A Constant Function

Let us take

$$f(x) = 1 = x^0 \tag{5.11}$$

The ordinary derivative is

$$Df(x) = 0 \tag{5.12}$$

In the half-derivative formula let us substitute

$$\beta = 0$$

Then

$$D^{1/2}1 = \frac{\Gamma(1)}{\Gamma\left(\frac{1}{2}\right)} x^{-1/2}$$

Since

$$\Gamma(1) = 1$$

and

$$\Gamma\left(\frac{1}{2}\right) = \sqrt{\pi},$$

we obtain

$$\boxed{D^{1/2}1 = \frac{1}{\sqrt{\pi}} x^{-1/2}} \tag{5.13}$$

Thus, for the constant function the three expressions to be compared are

$$f(x) = 1,$$

$$D^{1/2} f(x) = \frac{1}{\sqrt{\pi}} x^{-1/2},$$

$$Df(x) = 0 \tag{5.14}$$

Now let us check the second half-step. From the general formula

$$D^{1/2} x^{-1/2} = \frac{\Gamma\left(\frac{1}{2}\right)}{\Gamma(0)} x^{-1} \tag{5.15}$$

$\Gamma(0)$ is at a pole, while the reciprocal Gamma function $1/\Gamma(z)$ vanishes at $z = 0$. Therefore

$$\boxed{D^{1/2} x^{-1/2} = 0} \tag{5.16}$$

Accordingly,

$$D^{1/2}\left(D^{1/2} 1\right) = D^{1/2}\left(\frac{1}{\sqrt{\pi}} x^{-1/2}\right) \qquad = \frac{1}{\sqrt{\pi}} D^{1/2} x^{-1/2} \qquad = 0$$

Thus,

$$\boxed{D^{1/2}\left(D^{1/2} 1\right) = 0 = D1} \tag{5.17}$$

***Graphical behavior***

The function

$$\frac{1}{\sqrt{\pi}} x^{-1/2} = \frac{1}{\sqrt{\pi x}}$$

is positive on $x > 0$. Moreover,

when $x$ approaches zero from the right,

$$\frac{1}{\sqrt{\pi x}} \to +\infty.$$

while

when $x \to +\infty$,

$$\frac{1}{\sqrt{\pi x}} \to 0$$

Therefore the graph of the half-derivative of the constant function has vertical asymptotic behavior at $x = 0$ and gradually approaches zero as it moves to the right.

This example clearly shows the difference between the Riemann–Liouville-type half-derivative and the ordinary derivative: the first derivative of a constant is zero, but the first half-step is not zero.

## 5.4. Example II — The Function $f(x) = x$

Let us take

$$f(x) = x \tag{5.18}$$

The ordinary derivative is

$$Df(x) = 1 \tag{5.19}$$

For the half-derivative, in formula (5.2) let us substitute

$$\beta = 1$$

Then

$$D^{1/2}x = \frac{\Gamma(2)}{\Gamma\left(\frac{3}{2}\right)}x^{1/2}$$

We know that

$$\Gamma(2) = 1$$

and

$$\Gamma\left(\frac{3}{2}\right) = \frac{1}{2}\Gamma\left(\frac{1}{2}\right) = \frac{\sqrt{\pi}}{2}$$

Therefore

$$D^{1/2}x = \frac{1}{\sqrt{\pi}/2}x^{1/2},$$

that is

$$\boxed{D^{1/2}x = \frac{2}{\sqrt{\pi}}x^{1/2}} \tag{5.20}$$

Thus, we have three functions:

$$f(x) = x,$$

$$D^{1/2}f(x) = \frac{2}{\sqrt{\pi}}\sqrt{x},$$

$$Df(x) = 1 \tag{5.21}$$

***Verification of the second half-step***

We know that

$$D^{1/2}x^{1/2} = \frac{\sqrt{\pi}}{2}$$

Therefore

$$D^{1/2}\left(D^{1/2}x\right) = D^{1/2}\left(\frac{2}{\sqrt{\pi}}x^{1/2}\right) \qquad = \frac{2}{\sqrt{\pi}}D^{1/2}x^{1/2} \qquad = \frac{2}{\sqrt{\pi}} \cdot \frac{\sqrt{\pi}}{2} \qquad = 1$$

Thus,

$$\boxed{D^{1/2}\left(D^{1/2}x\right) = 1 = Dx} \tag{5.22}$$

***Graphical behavior***

The initial function $x$ is linear. Its half-derivative

$$\frac{2}{\sqrt{\pi}}\sqrt{x}$$

increases, but its rate of increase gradually decreases. The ordinary derivative, however, is a constant function:

$$Df(x) = 1$$

These three graphs differ from one another. The half-derivative's graph is neither tangent to $f(x) = x$ nor the graph of its ordinary derivative.

## 5.5. Example III — The Function $f(x) = x^2$

Let us take

$$f(x) = x^2 \tag{5.23}$$

Its ordinary derivative is

$$Df(x) = 2x \tag{5.24}$$

For the half-derivative, in formula (5.2) let us substitute

$$\beta = 2$$

Then

$$D^{1/2}x^2 = \frac{\Gamma(3)}{\Gamma\left(\frac{5}{2}\right)}x^{3/2}$$

We know that

$$\Gamma(3) = 2$$

and

$$\Gamma\left(\frac{5}{2}\right) = \frac{3}{2}\Gamma\left(\frac{3}{2}\right) = \frac{3}{2} \cdot \frac{\sqrt{\pi}}{2} = \frac{3\sqrt{\pi}}{4}$$

Therefore

$$D^{1/2}x^2 = \frac{2}{3\sqrt{\pi}/4}x^{3/2},$$

from which

$$\boxed{D^{1/2}x^2 = \frac{8}{3\sqrt{\pi}}x^{3/2}} \tag{5.25}$$

Thus, the functions to be compared are

$$f(x) = x^2,$$

$$D^{1/2}f(x) = \frac{8}{3\sqrt{\pi}}x^{3/2},$$

$$Df(x) = 2x \tag{5.26}$$

***Verification of the second half-step***

We know that

$$D^{1/2}x^{3/2} = \frac{3\sqrt{\pi}}{4}x$$

Therefore

$$D^{1/2}\left(D^{1/2}x^2\right) = D^{1/2}\left(\frac{8}{3\sqrt{\pi}}x^{3/2}\right) \qquad = \frac{8}{3\sqrt{\pi}}D^{1/2}x^{3/2} \qquad = \frac{8}{3\sqrt{\pi}} \cdot \frac{3\sqrt{\pi}}{4}x \qquad = 2x$$

Thus,

$$\boxed{D^{1/2}\left(D^{1/2}x^2\right) = 2x = Dx^2} \tag{5.27}$$

### 5.6. The Transformation of $x^2$ for a General $\alpha$-Order

Apart from the half-order, it is useful to see how the result of acting on $x^2$ changes as $\alpha$ gradually varies from 0 to 1.

By the general formula

$$D^{\alpha} x^2 = \frac{\Gamma(3)}{\Gamma(3-\alpha)} x^{2-\alpha}$$

Since

$$\Gamma(3) = 2,$$

we obtain

$$\boxed{D^{\alpha} x^2 = \frac{2}{\Gamma(3-\alpha)} x^{2-\alpha}} \tag{5.28}$$

Table 5: The action of $D^{\alpha}$ on $x^2$ for selected orders

| $\alpha$ | $D^{\alpha} x^2$ | **Degree** |
|---|---|---|
| 0 | $x^2$ | 2 |
| $\frac{1}{4}$ | $\frac{2}{\Gamma\left(\frac{11}{4}\right)} x^{7/4}$ | $\frac{7}{4}$ |
| $\frac{1}{2}$ | $\frac{8}{3\sqrt{\pi}} x^{3/2}$ | $\frac{3}{2}$ |
| $\frac{3}{4}$ | $\frac{2}{\Gamma\left(\frac{9}{4}\right)} x^{5/4}$ | $\frac{5}{4}$ |
| 1 | $2x$ | 1 |

When

$$\alpha = 0,$$

we obtain the initial function:

$$D^0 x^2 = x^2$$

When

$$\alpha = 1,$$

we obtain the ordinary derivative:

$$D^1 x^2 = 2x$$

Intermediate values

$$0 < \alpha < 1$$

give the corresponding intermediate-order operator transformations.

For a fixed $x > 0$ and for those $\alpha$ for which the Gamma-coefficient is defined, the expression

$$\alpha \mapsto D^{\alpha} x^2$$

varies continuously. However this fact should not be equated with a complete theorem of continuity of $D^{\alpha}$-operators on a general functional space. Here we are talking only about the specific monomial formula.

## 5.7. Example IV — The Polynomial P(x)=1+3x+2x$^2$

Let us consider the polynomial

$$P(x) = 1 + 3x + 2x^2 \tag{5.29}$$

Its ordinary derivative is

$$DP(x) = 3 + 4x \tag{5.30}$$

Let us find the half-derivative term-by-term:

$$D^{1/2} P(x) = D^{1/2} 1 + 3 D^{1/2} x + 2 D^{1/2} x^2$$

We already know:

$$D^{1/2} 1 = \frac{1}{\sqrt{\pi}} x^{-1/2},$$

$$D^{1/2} x = \frac{2}{\sqrt{\pi}} x^{1/2},$$

and

$$D^{1/2} x^2 = \frac{8}{3\sqrt{\pi}} x^{3/2}$$

Therefore

$$D^{1/2} P(x) = \frac{1}{\sqrt{\pi}} x^{-1/2} + 3 \cdot \frac{2}{\sqrt{\pi}} x^{1/2} + 2 \cdot \frac{8}{3\sqrt{\pi}} x^{3/2} \qquad = \frac{1}{\sqrt{\pi}} x^{-1/2} + \frac{6}{\sqrt{\pi}} x^{1/2} + \frac{16}{3\sqrt{\pi}} x^{3/2}$$

Thus,

$$\boxed{D^{1/2} P(x) = \frac{1}{\sqrt{\pi}} \left( x^{-1/2} + 6 x^{1/2} + \frac{16}{3} x^{3/2} \right)} \tag{5.31}$$

After the second half-step we obtain

$$D^{1/2} \left( D^{1/2} P \right) = 0 + 3 + 4x \qquad = 3 + 4x$$

Therefore

$$\boxed{D^{1/2}\left(D^{1/2}P\right) = DP} \tag{5.32}$$

Because of the presence of the constant term in the polynomial, $D^{1/2}P\,(x)$ contains the term

$$x^{-1/2}$$

and is singular at $x = 0$.

If instead we take the polynomial

$$Q\,(x) = 3x + 2x^2,$$

for which

$$Q\,(0) = 0,$$

then

$$D^{1/2}Q\,(x) = \frac{6}{\sqrt{\pi}}x^{1/2} + \frac{16}{3\sqrt{\pi}}x^{3/2}, \tag{5.33}$$

and this function tends to zero as $x$ approaches zero from the right.

Thus, in visualization the presence or absence of the constant term significantly affects the behavior of the graph near zero.

## 5.8. A Summary Table of Formulas

Table 6: Initial Function, First Half-Step, and Final Derivative

| $f\,(x)$ | $D^{1/2}f\,(x)$ | $D^{1/2}\left(D^{1/2}f\,(x)\right)$ | $Df\,(x)$ |
|---|---|---|---|
| $1$ | $\frac{1}{\sqrt{\pi}}x^{-1/2}$ | $0$ | $0$ |
| $x$ | $\frac{2}{\sqrt{\pi}}x^{1/2}$ | $1$ | $1$ |
| $x^2$ | $\frac{8}{3\sqrt{\pi}}x^{3/2}$ | $2x$ | $2x$ |
| $x^3$ | $\frac{16}{5\sqrt{\pi}}x^{5/2}$ | $3x^2$ | $3x^2$ |
| $1 + 3x + 2x^2$ | $\frac{1}{\sqrt{\pi}}\left(x^{-1/2} + 6x^{1/2} + \frac{16}{3}x^{3/2}\right)$ | $3 + 4x$ | $3 + 4x$ |

The table demonstrates two distinct facts:

- the first half-step result is, in general, neither the initial function nor its ordinary derivative;
- after the second half-step the ordinary first derivative is obtained exactly.

## 5.9. Numerical Verification for $f\,(x) = x^2$

The half-derivative of the function

$$f(x) = x^2$$

is

$$D^{1/2} f(x) = \frac{8}{3\sqrt{\pi}} x^{3/2}$$

The approximate value of the coefficient is

$$\frac{8}{3\sqrt{\pi}} \approx 1.504505556 \tag{5.34}$$

Therefore

$$D^{1/2} x^2 \approx 1.504505556\, x^{3/2}$$

Table 7: Numerical Values

| $x$ | $f(x) = x^2$ | $D^{1/2} f(x)$ | $D^{1/2}\left(D^{1/2} f(x)\right)$ | $Df(x) = 2x$ |
|---|---|---|---|---|
| 0.25 | 0.062500 | 0.188063 | 0.500000 | 0.500000 |
| 0.50 | 0.250000 | 0.531923 | 1.000000 | 1.000000 |
| 1.00 | 1.000000 | 1.504506 | 2.000000 | 2.000000 |
| 2.00 | 4.000000 | 4.255384 | 4.000000 | 4.000000 |
| 4.00 | 16.000000 | 12.036044 | 8.000000 | 8.000000 |

The table shows once again that the value of $D^{1/2} f(x)$ does not always lie numerically between $f(x)$ and $Df(x)$. For example, at $x = 2$ we have

$$f(2) = 4, \qquad D^{1/2} f(2) \approx 4.255384, \qquad Df(2) = 4$$

Thus, "intermediate" means the order of the operator and not the necessary numerical arrangement of the graphs.

The values of the second half-step, however, coincide with $2x$ at all the points indicated. Theoretically the difference is exactly zero; in a computer calculation there may remain only an insignificant error caused by rounding.

## 5.10. Graphical Interpretation

For visual teaching it is useful to plot three graphs in same coordinate system:

$$y = f(x),$$

$$y = D^{1/2} f(x),$$

$$y = Df(x)$$

For example, in the case $f(x) = x^2$, the graphs to be compared are

$$y = x^2,$$

$$y = \frac{8}{3\sqrt{\pi}} x^{3/2},$$

$$y = 2x \tag{5.35}$$

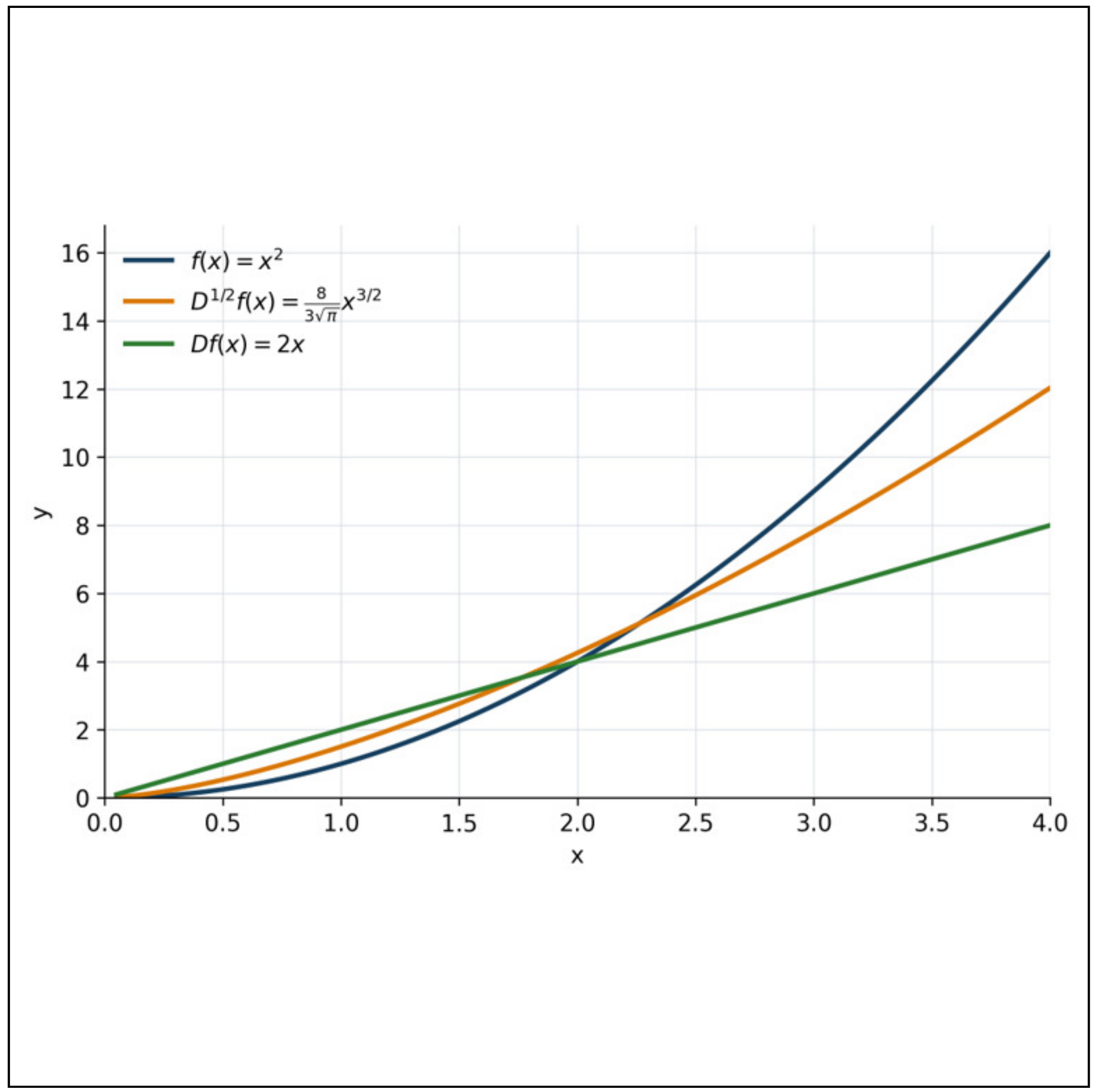

Figure 1: Comparison of $f(x) = x^2$, $D^{1/2}f(x)$, and $Df(x)$ for $x > 0$.

For the same example, separately one could plot

$$y = D^{1/2}\left(D^{1/2} f\left(x\right)\right)$$

and

$$y = Df\left(x\right)$$

These two graphs must coincide with each other.

In visualization it is important to observe the following distinctions.

#### *1. The graph of the half-derivative is not a tangent*

$D^{1/2}f(x)$ is a new function. Its graph does not represent the tangent line drawn at any point of $f(x)$.

#### *2. The half-derivative is not the ordinary derivative*

In general,

$$D^{1/2}f(x) \neq Df(x)$$

For the admissible functions considered in this chapter, equality is obtained only after the second half-step action:

$$D^{1/2}\left(D^{1/2}f\right) = Df$$

#### *3. The graph is not required to lie between the other two graphs*

The intermediacy of the operator order does not imply pointwise numerical intermediacy.

#### *4. In physical models the dimensions must be taken into account*

If $f$, $D^{1/2}f$ and $Df$ describe physical quantities, they may have different physical dimensions. In such a case, before direct comparison on the same axes, an appropriate normalization or transition to dimensionless variables is needed.

## 5.11. The Domain of Definition and the Question of $x = 0$

The general monomial formula contains the power

$$x^{\beta-\alpha}$$

Therefore behavior near $x = 0$ depends on the difference

$$\beta - \alpha$$

If

$$\beta - \alpha > 0,$$

then

as $x$ approaches zero from the right,

$$x^{\beta-\alpha} \to 0$$

If

$$\beta - \alpha = 0,$$

the power part of the result is constant.

If

$$\beta - \alpha < 0,$$

then $x^{\beta-\alpha}$ is, in general, singular at $x = 0$.

An exception can arise when the Gamma-function coefficient vanishes because of the pole structure of the reciprocal Gamma function.

For example,

$$D^{1/2}1 = \frac{1}{\sqrt{\pi}}x^{-1/2}$$

is not defined at $x = 0$. Therefore in constructing its graph one should not use an interval that begins exactly at zero.

We could take, for example,

$$x \in [0.01, 5]. \tag{5.36}$$

This issue is methodologically significant: after the action of the fractional operator the function's domain of definition and its behavior near zero may change.

## 5.12. A Computational Algorithm

For numerical computation of the fractional derivative of a monomial the following data are needed:

$$\alpha, \qquad \beta, \qquad x > 0$$

The algorithm can be formulated in the following steps:

- compute the coefficient

$$C_\alpha(\beta) = \frac{\Gamma(\beta+1)}{\Gamma(\beta+1-\alpha)};$$

- compute the new power

$$\beta - \alpha;$$

- obtain

$$D^{\alpha}x^{\beta} = C_{\alpha}(\beta)\,x^{\beta-\alpha}$$

For a polynomial

$$P(x) = \sum_{j=0}^{N} a_j x^j$$

the same action is performed for each term:

$$D^{\alpha}P(x) = \sum_{j=0}^{N} a_j C_{\alpha}(j)\,x^{j-\alpha}$$

In a computer program it is desirable to use not only $\Gamma(z)$, but also the reciprocal Gamma function

$$\mathrm{rgamma}(z) = \frac{1}{\Gamma(z)}$$

In this case the coefficient can be written in the form

$$C_{\alpha}(\beta) = \Gamma(\beta+1)\,\mathrm{rgamma}(\beta+1-\alpha) \tag{5.37}$$

This notation is particularly convenient when the denominator's Gamma function falls at a pole, since *rgamma* vanishes at the corresponding points.

## 5.13. Python Code for Visualization and Verification

The code below plots

$$f(x) = x^2,$$

its half-derivative, and its ordinary derivative. At the same time it numerically checks the second half-step action.

```
import numpy as np
import matplotlib.pyplot as plt
from scipy.special import gamma, rgamma


def d_fractional(x, beta, alpha):
    """Compute D^alpha(x^beta) for x > 0."""
    x = np.asarray(x, dtype=float)
    if np.any(x <= 0):
        raise ValueError("All x-values must be positive.")

    coefficient = gamma(beta + 1.0) * rgamma(beta + 1.0 - alpha)
    return coefficient * x ** (beta - alpha)


beta, alpha = 2.0, 0.5
x = np.linspace(0.01, 5.0, 500)
```

```
f = x ** beta
d_half = d_fractional(x, beta, alpha)
d_first = beta * x ** (beta - 1.0)

c_half = gamma(beta + 1.0) * rgamma(beta + 1.0 - alpha)
d_half_twice = c_half * d_fractional(x, beta - alpha, alpha)

error = np.max(np.abs(d_half_twice - d_first))
print("Maximum composition error:", error)

plt.figure()
plt.plot(x, f, label=r"$f(x)=x^2$")
plt.plot(x, d_half, label=r"$D^{1/2}f(x)$")
plt.plot(x, d_first, label=r"$Df(x)=2x$")
plt.xlabel("x")
plt.ylabel("Function value")
plt.legend()
plt.grid(True)
plt.tight_layout()
plt.show()
```

Theoretically the program should confirm

$$D^{1/2}\left(D^{1/2}x^2\right) = 2x$$

In the computer computation the maximum error obtained should be very close to zero and caused by rounding of floating-point numbers.

Program verification does not replace mathematical proof. It only confirms that the formula has been correctly transferred into the computational algorithm.

## 5.14. The Possibility of Interactive Visualization

In an interactive environment, for example using Streamlit, it is possible for the user to change:

$$\alpha, \qquad \beta, \qquad x_{min}, \qquad x_{max}$$

The program could simultaneously show:

$$f\left(x\right) = x^{\beta},$$

$$D^{\alpha} f\left(x\right) = \frac{\Gamma\left(\beta+1\right)}{\Gamma\left(\beta+1-\alpha\right)} x^{\beta-\alpha},$$

and

$$Df\left(x\right) = \beta x^{\beta-1}$$

If

$$0 \leq \alpha \leq 1,$$

by moving the slider the student will see how the monomial formula changes:

from the form

$$D^0 x^\beta = x^\beta$$

to the form

$$D^1 x^\beta = \beta x^{\beta-1}$$

The program could be supplemented with:

- input of the coefficients of a polynomial;
- comparison of the Riemann–Liouville and Caputo results;
- automatic computation of the second half-step;
- display of the composition error;
- a warning connected with the function's domain of definition;
- control of the pole cases of the Gamma function.

The purpose of the interactive model is a visual presentation of the formulas obtained on the monomial space. Such a program does not by itself prove continuity of the operator on a general functional space, nor its non-local or physical properties.

### 5.15. Methodological Significance and the Limits of Computational Verification

The examples discussed in this chapter demonstrate several important points.

First, the fractional derivative represents a new function obtained from a function:

$$f \mapsto D^\alpha f$$

Second, for the admissible monomial and polynomial functions studied here, the first half-order action differs from the ordinary derivative, but the second action is compositionally linked with it:

$$D^{1/2}\left(D^{1/2} f\right) = Df$$

Third, the action on polynomial initial functions can be computed term-by-term.

Fourth, the constant term, under the Riemann–Liouville-type action, produces the singular term

$$\frac{a_0}{\sqrt{\pi}} x^{-1/2}$$

Fifth, graphical representation helps the student distinguish:

- the initial function;
- the fractional-order transformation;
- the ordinary derivative.

Sixth: numerical verification is useful for confirming the formulas, the program code, and specific examples, but it does not represent a proof of a general theorem.

In particular, graphs and numerical examples do not prove:

- the existence of an integral or convolution kernel;
- the non-locality of the operator;
- the property of causality;
- continuity on a broad functional space;
- solvability of initial-boundary value problems;
- applicability in a physical model.

These issues go beyond the bounds of the present monomial and polynomial construction.

## 5.16. Summary of the Chapter

In this chapter the formulas obtained in the previous chapters were presented with concrete computational and visual examples.

For the constant function we obtained

$$D^{1/2}1 = \frac{1}{\sqrt{\pi}}x^{-1/2},$$

and

$$D^{1/2}\left(D^{1/2}1\right) = 0$$

For the function $f(x) = x$ we obtained

$$D^{1/2}x = \frac{2}{\sqrt{\pi}}x^{1/2},$$

and

$$D^{1/2}\left(D^{1/2}x\right) = 1$$

For the function $f(x) = x^2$ we obtained

$$D^{1/2}x^2 = \frac{8}{3\sqrt{\pi}}x^{3/2},$$

and

$$D^{1/2}\left(D^{1/2}x^2\right) = 2x$$

For the polynomial

$$P(x) = 1 + 3x + 2x^2$$

we obtained

$$D^{1/2}P\left(x\right)=\frac{1}{\sqrt{\pi}}\left(x^{-1/2}+6x^{1/2}+\frac{16}{3}x^{3/2}\right),$$

while by the second half-step action

$$D^{1/2}\left(D^{1/2}P\right)=3+4x=DP$$

The numerical table and the computer algorithm showed how the formulas obtained can be practically verified. Graphical discussion clarified that the fractional derivative is not the slope of a tangent, and that its "intermediacy" concerns the order of the operator, not the necessary position of the graph.

In the next chapter the algebraic, operator, and methodological results of the whole construction will be summarized, the limits of the results obtained will be clearly delineated, and directions for further research will be formulated.

# Chapter VI — Methodological Summary, the Limits of the Construction, and Directions for Further Development

## 6.1. The Overall Picture of the Construction and the Logical Chain

In the present paper the fractional-order operator on monomials was built step by step, without the prior use of an integral definition or a limiting process. Every compositional step in the following logical chain is understood on the corresponding admissible monomial domain:

**The ordinary derivative and the restricted monomial-level composition requirement** $D^{1/2}\circ D^{1/2}=D$ [26] **→ coefficient recurrence** $c\left(\beta\right)c\left(\beta-\frac{1}{2}\right)=\beta$ **[26] → one-step recurrence on integer arguments (1.1) → functional equation and its canonical Gamma-function solution** $c\left(\beta\right)=\Gamma\left(\beta+1\right)/\Gamma\left(\beta+\frac{1}{2}\right)$ **(1.8, canonicity — Section 1.10) → division into** $m$ **equal steps and rational orders (3.1.2–3.1.3) → the general coefficient condition following from the operator addition law** $C_{\gamma}\left(\beta\right)C_{\alpha}\left(\beta-\gamma\right)=C_{\alpha+\gamma}\left(\beta\right)$ **(3.22) → the general-order monomial formula**

$$\boxed{D^{\alpha}x^{\beta}=\frac{\Gamma\left(\beta+1\right)}{\Gamma\left(\beta+1-\alpha\right)}x^{\beta-\alpha}},x>0 \tag{3.37}$$

**→ linear extension to polynomials and comparison with the Riemann–Liouville and Caputo formulas (Chapter IV).**

Two structural observations in this chain deserve emphasis. First, the two discrete coefficient families of integer and half-integer powers established in Part I [26] are special cases of one Gamma-function formula, where the factor $\sqrt{\pi}$ appears in the numerator in one chain and in the denominator in the other, and cancels under the double composition — this is a structural, and not an accidental numerical, coincidence. Second, the double-factorial forms established in Part I [26] represent a bridge between the elementary recurrence and the Gamma-function form — the Gamma function here is not used as the initial definition, but appears as a compact functional notation of an already obtained discrete structure.

## 6.2. Classification of Results by Order

Formula (3.37) uniformly covers three cases (detailed derivation — Sections 3.11 and 3.16): at nonnegative integer $\alpha = m$ it reverts to ordinary repeated differentiation (3.39); at positive non-integer $\alpha$ it coincides with the left-sided Riemann–Liouville monomial formula with lower limit 0; while at negative $\alpha = -\rho$ it gives the monomial form of the corresponding-order Riemann–Liouville-type integral (3.45), which was also independently confirmed by means of the Beta function (3.51).

## 6.3. Extension to Polynomials and Comparison with the Classical Operators

The operator obtained on monomials was linearly extended to polynomial initial functions (4.11), (4.20). The decisive difference between the Riemann–Liouville and Caputo [20–24] forms lies in the fate of the constant and low-degree terms: the present construction gives a nonzero result on the constant as well (4.17)–(4.19), and thereby follows the Riemann–Liouville rule, while agreement with the Caputo form is obtained only for sufficiently high-degree terms (4.21)–(4.24).

Table 8: The Relation of the Obtained Formula to the Classical Rules

| **Case** | **Obtained algebraic–operator result** | **Riemann–Liouville form** | **Caputo form** |
|---|---|---|---|
| $x^\beta$, admissible $\beta$ | $\frac{\Gamma(\beta+1)}{\Gamma(\beta+1-\alpha)}x^{\beta-\alpha}$ | coincides | depends on the conditions |
| constant 1 | $\frac{x^{-\alpha}}{\Gamma(1-\alpha)}$ | coincides | 0 |
| $x^j$, $j \geq m$ | Gamma-function formula | coincides | coincides |
| $x^j$, $0 \leq j < m$ | Gamma-function formula | coincides | 0 |

## 6.4. Visual and Technological Aspects

The computational examples of Chapter V (see 5.2–5.10) showed that the simultaneous graphical comparison of the three functions — $f$, $D^{1/2}f$, and $Df$ — clearly demonstrates that the graph of the half-derivative is neither the initial function, nor the ordinary derivative, nor a tangent line, but the graph of an independent function which, under a double action, passes exactly into $Df$. In the interactive program (`https://kapanadze-theory.streamlit.app`), by changing the parameters $\alpha$, $\beta$, $x_{min}$, $x_{max}$ it is possible to observe this dependence in real time. Visualization remains a means of explaining and verifying the formulas, and not a substitute for theoretical proof.

## 6.5. The Exact Bounds of the Results Obtained

In the present paper the following results have been obtained.

### *Proven*

- From the half-step composition requirement the coefficient recurrence is obtained

$$c\,(\beta)\,c\left(\beta - \frac{1}{2}\right) = \beta$$

- On nonnegative integer and half-integer powers two interlinked coefficient families are constructed.

- Their canonical functional continuation is expressed by the formula

$$c(\beta) = \frac{\Gamma(\beta+1)}{\Gamma\left(\beta+\frac{1}{2}\right)}$$

- From the constructive motivation of rational orders and the operator addition law the general coefficient is obtained

$$C_\alpha(\beta) = \frac{\Gamma(\beta+1)}{\Gamma(\beta+1-\alpha)}$$

- On monomials, on the corresponding domain, the following holds

$$D^\alpha D^\gamma = D^{\alpha+\gamma}$$

- For integer orders the formula reverts to ordinary repeated differentiation.
- On positive non-integer orders the formula coincides with the left-sided Riemann–Liouville monomial rule with lower limit 0.
- The operator extends linearly to polynomial initial functions.

***Not proven***

The present paper does not construct:

- the Riemann–Liouville integral kernel;
- the complete Caputo operator construction;
- an independent algebraic mechanism of non-locality;
- a convolution representation;
- a general rule for the Laplace or Fourier transform;
- the complete domain of the operator on a broad functional space;
- a functional-analytic theory of continuity, closedness, or generation;
- a theory of existence and uniqueness for fractional differential equations;
- direct applicability in physical or biological models.

This delimitation is of principled importance: coincidence with the classical formula at the monomial level does not mean the automatic reconstruction of the complete theory of the corresponding classical operator.

## 6.6. Methodological Significance

The main methodological value of the present approach is the step-by-step construction of the result.

In classical teaching the formula for the fractional derivative of a monomial is often written directly:

$$D^{\alpha}x^{\beta} = \frac{\Gamma(\beta+1)}{\Gamma(\beta+1-\alpha)}x^{\beta-\alpha}$$

In this case the student may be left with natural questions:

- why does the degree decrease exactly by $\alpha$?
- why does the ratio of the Gamma functions appear?
- why do the orders of the operators add?
- how is the half-derivative connected to the ordinary derivative?

The present construction answers these questions in the following sequence:

**Ordinary derivative → half-step requirement → recurrence → two coefficient chains → double factorials → functional continuation → Gamma function → general order.**

Thus, the student obtains the final formula not as a pre-given rule, but as the result of consistent operator and algebraic requirements.

## 6.7. Directions for Further Development

The monomial construction obtained suggests several directions for further development.

#### *1. Transcendental Functions*

An important question is the extension of the formula to functions such as

$$e^{x}, \qquad lnx, \qquad sinx, \qquad cosx$$

For this, series, functional equations, or corresponding operator properties could be used. Particularly important will be the justification of convergence and of term-by-term action.

#### *2. A Caputo-Type Modification*

The difference between Riemann–Liouville and Caputo lies in the low-degree terms. Therefore it could be separately investigated whether, on the monomial space, a modified operator can be defined which annihilates the terms

$$1, x, \ldots, x^{m-1}$$

and preserves the same Gamma-coefficient for higher degrees.

#### *3. Extension of the Functional Space*

It is necessary to define such a space $\mathcal{V}$ on which

$$D^{\alpha} : \mathcal{V} \to \mathcal{V}$$

or the action between corresponding spaces will be strictly formulated.

#### *4. The Question of Non-Locality*

One of the basic properties of the classical fractional derivatives is non-locality. In the monomial coefficient formula this property is not directly visible. Therefore it is necessary to investigate whether an integral, convolution, or other non-local representation can be obtained from the algebraic–operator construction.

#### *5. The Theory of Transforms*

A further task is establishing the connection with the Laplace and Fourier transforms. For example, it should be investigated under what conditions a formula of the form

$$\mathcal{L}\left(D^{\alpha} f\right)(s) = s^{\alpha} F(s) - S_0(s) \qquad (6.50)$$

can be obtained. Here $S_0(s)$ denotes the terms determined by the initial data.

#### *6. Fractional Differential Equations*

After a broader definition of the operator it will become possible to consider such equations as

$$D^{\alpha} y(t) = \lambda y(t), \qquad (6.1)$$

or

$$D_t^{\alpha} u(x,t) = K \frac{\partial^2 u}{\partial x^2} \qquad (6.2)$$

At this stage such an application does not yet follow from the monomial formula alone.

## 6.8. Final Conclusion

The purpose of the present paper has not been the discovery of a new classical formula of fractional differentiation. The purpose was to obtain the known monomial Gamma-formula by a different route — without the prior use of the integral definition, on the basis of algebraic recurrence and operator composition.

In Part I the half-order operator was built, on the initial admissible monomial class, from the requirement

$$D^{1/2} \circ D^{1/2} = D,$$

*This relation is used throughout the paper only in the restricted sense stated above; it is not claimed as a universal semigroup identity on arbitrary functions.*

which gave rise to a coefficient recurrence, two interlinked chains, double factorials, and a normalization constant fixed by the compatibility condition.

In the present paper, this structure was extended to a functional equation, from which the canonical form of the half-derivative was obtained

$$D^{1/2}x^{\beta} = \frac{\Gamma(\beta+1)}{\Gamma\left(\beta+\frac{1}{2}\right)}x^{\beta-\frac{1}{2}} \tag{6.3}$$

Then from the operator addition law the general-order formula was obtained

$$\boxed{D^{\alpha}x^{\beta} = \frac{\Gamma(\beta+1)}{\Gamma(\beta+1-\alpha)}x^{\beta-\alpha}} \tag{6.4}$$

This result, on monomials:

- contains the ordinary derivative as the special case of integer order;
- contains the half-derivative as the case $\alpha = \frac{1}{2}$;
- coincides with the left-sided Riemann–Liouville monomial formula on the corresponding domain;
- extends linearly to polynomial initial functions;
- preserves the addition law of operator orders on the corresponding domain.

At the same time the paper clearly bounds the scale of the result obtained: what is constructed is an algebraic–operator family acting on the system of monomials, and not the complete integral and functional-analytic theory of fractional calculus.

The main methodological idea of the construction obtained can be summarized as follows:

**The classical monomial formula of fractional order is obtained as the result of operator composition and algebraic continuation.**

Thus, the Gamma function does not appear in the final formula as an opaque, pre-given instrument. It emerges as the natural and canonical functional form of an already discovered coefficient regularity.